\documentclass{amsart}
\RequirePackage[OT1]{fontenc}
\RequirePackage{amsthm,amsmath,amssymb}

\usepackage{epsfig,graphicx,color}
\definecolor{darkblue}{rgb}{.2, 0.2,.8}
\definecolor{darkgreen}{rgb}{0,0.5,0.3}
\definecolor{darkred}{rgb}{.8, .1,.1}

\newcommand{\asy}{asymptotic}

\newcommand{\ts}{time series}
\newcommand{\tsa}{\ts\ analysis}

\newtheorem{lemma}{Lemma}[section]

\newtheorem{theorem}[lemma]{Theorem}

\newtheorem{proposition}[lemma]{Proposition}
\newtheorem{definition}[lemma]{Definition}
\newtheorem{corollary}[lemma]{Corollary}
\newtheorem{example}[lemma]{Example}
\newtheorem{exercise}[lemma]{Exercise}
\newtheorem{remark}[lemma]{Remark}
\newtheorem{fig}[lemma]{Figure}
\newtheorem{tab}[lemma]{Table}

\newcommand{\MC}{Markov chain}

\newcommand{\bfR}{{\bf R}}
\newcommand{\bfC}{{\bf C}}

\newcommand{\RV}{{\rm RV}}

\newcommand{\bth}{\begin{theorem}}
\newcommand{\ethe}{\end{theorem}}
\newcommand{\sv}{stochastic volatility}
\newcommand{\bre}{\begin{remark}\em }
\newcommand{\ere}{\end{remark}}

\newcommand{\ble}{\begin{lemma}}
\newcommand{\ele}{\end{lemma}}
\newcommand{\sre}{stochastic recurrence equation}
\newcommand{\pp}{point process}
\newcommand{\bde}{\begin{definition}}
\newcommand{\ede}{\end{definition}}
\newcommand{\bco}{\begin{corollary}}
\newcommand{\eco}{\end{corollary}}

\newcommand{\bpr}{\begin{proposition}}
\newcommand{\epr}{\end{proposition}}

\newcommand{\bexer}{\begin{exercise}}
\newcommand{\eexer}{\end{exercise}}

\newcommand{\bexam}{\begin{example}}
\newcommand{\eexam}{\end{example}}

\newcommand{\bfi}{\begin{fig}}
\newcommand{\efi}{\end{fig}}

\newcommand{\btab}{\begin{tab}}
\newcommand{\etab}{\end{tab}}

\newcommand{\rv}{random variable}

\newcommand{\bfTh}{\mbox{\boldmath$\Theta$}}

\newcommand{\rhs}{right-hand side}

\newcommand{\dint}{\displaystyle\int}

\newcommand{\dprod}{\displaystyle\prod}

\newcommand{\beao}{\begin{eqnarray*}}
\newcommand{\eeao}{\end{eqnarray*}\noindent}

\newcommand{\beam}{\begin{eqnarray}}
\newcommand{\eeam}{\end{eqnarray}\noindent}

\newcommand{\beqq}{\begin{equation}}
\newcommand{\eeqq}{\end{equation}\noindent}

\newcommand{\bce}{\begin{center}}
\newcommand{\ece}{\end{center}}

\newcommand{\barr}{\begin{array}}
\newcommand{\earr}{\end{array}}

\newcommand{\std}{\stackrel{d}{\rightarrow}}

\newcommand{\stv}{\stackrel{v}{\rightarrow}}
\newcommand{\stw}{\stackrel{w}{\rightarrow}}

\newcommand{\eqd}{\stackrel{d}{=}}

\newcommand{\vague}{\stackrel{\lower0.2ex\hbox{$\scriptscriptstyle
                    \it{v} $}}{\rightarrow}}
\newcommand{\weak}{\stackrel{\lower0.2ex\hbox{$\scriptscriptstyle
                    \it{w} $}}{\rightarrow}}
\newcommand{\what}{\stackrel{\lower0.2ex\hbox{$\scriptscriptstyle
                    \it{\hat{w}} $}}{\rightarrow}}

\newcommand{\bdis}{\begin{displaymath}}
\newcommand{\edis}{\end{displaymath}\noindent}

\newcommand{\Prob}{\operatorname{P}}

\newcommand{\R}{\mathbb{R}}
\newcommand{\Ro}{\mathbb{R}_{\mathbf{0}}}

\newcommand{\nto}{n\to\infty}
\newcommand{\kto}{k\to\infty}
\newcommand{\xto}{x\to\infty}

\newcommand{\wt}{\widetilde}
\newcommand{\wh}{\widehat}
\newcommand{\vep}{\varepsilon}

\newcommand{\regvary}{regularly varying}
\newcommand{\slvary}{slowly varying}
\newcommand{\regvar}{regular variation}

\newcommand{\bbr}{{\mathbb R}}

\newcommand{\bbz}{{\mathbb Z}}

\newcommand{\bbs}{{\mathbb S}}

\newcommand{\con}{convergence}

\newcommand{\wrt}{with respect to}

\newcommand{\fct}{function}

\newcommand{\ds}{distribution}

\newcommand{\rep}{representation}

\newcommand{\seq}{sequence}

\newcommand{\pro}{probabilit}

\newcommand{\ms}{measure}

\newcommand{\bfx}{{\bf x}}
\newcommand{\bfX}{{\bf X}}
\newcommand{\bfB}{{\bf B}}
\newcommand{\bfY}{{\bf Y}}
\newcommand{\bfy}{{\bf y}}
\newcommand{\bfA}{{\bf A}}

\newcommand{\bfZ}{{\bf Z}}
\newcommand{\bfz}{{\bf z}}

\newcommand{\bfa}{{\bf a}}

\newcommand{\E }{\mathbb E}
\renewcommand{\P }{\mathbb P}

\newcommand{\1}{{\mathbf 1}}

\usepackage{bbm}

\def\Erw{\mathbb{E}}

\def\M{\mathbb{M}}

\def\Prob{\mathbb{P}} 

\def\R{\mathbb{R}}

\newcommand{\bX}{{\mathbf X}}
\newcommand{\bY}{{\mathbf Y}}
\newcommand{\bR}{{\mathbf R}}

\def\bs{{\bf s}}
\def\br{{\bf r}}
\def\bx{{\bf x}}
\def\by{{\bf y}}

\def\bz{{\bf z}}

\def\ba{{\bf a}}
\def\bb{{\bf b}}

\def\3{{\ss}}

\def\vag{\stackrel{v}{\to}}

\def\wh{\widehat}

\usepackage[top=3cm, bottom=3cm, outer=3.5cm, inner=3cm]{geometry}

\usepackage{soul,xcolor}

\def\diag{{\rm diag}}
\begin{document}

\title{Homogeneous mappings of regularly varying vectors}


\author{Piotr Dyszewski}

\address{Institute of Mathematics, University of Wroc\l{}aw\\ pl. Grunwaldzki 2/4,\\ 50-384 Wroc\l{}aw,\\ Poland}

\email{piotr.dyszewski@math.uni.wroc.pl}

\thanks{Piotr Dyszewski was partially supported by the National Science Centre, Poland
 (Sonata Bis, grant number DEC-2014/14/E/ST1/00588). This work was initiated while the first author was visiting the Department of Mathematics, University of Copenhagen in February 2018. He gratefully acknowledges financial support and hospitality.
}


%

\author{Thomas Mikosch}
\address{Department of Mathematics, University of Copenhagen\\ Universitetsparken 5,\\ DK-2100 Copenhagen,\\ Denmark}
\email{mikosch@math.ku.dk}
\thanks{Thomas Mikosch
is partially supported by an Alexander von Humboldt Research Award. He takes 
pleasure in thanking his colleagues at Mathematische Fakult\"at of 
Ruhruniversit\"at Bochum for hosting him December 2018 - May 2019.}



\date{\today}



\begin{abstract}
It is well known that the product of two independent
\regvary\ \rv s with the same tail index is again \regvary\ with this
index. In this paper, we provide sharp sufficient conditions for 
the \regvar\ property of product-type \fct s of \regvary\ random vectors, generalizing and extending the univariate theory in various directions. 
The main result is then applied to characterize the \regvar\ property of 
products of iid \regvary\ quadratic
random matrices and of solutions to affine \sre s 
under non-standard conditions.
\end{abstract}

\maketitle

\section{Introduction}\label{sec:introduction}

\subsection{Closure of \regvar\ under multiplication -- the univariate case}
	Consider a non-negative \rv\ $X$ and assume that $X$ is
	{\em \regvary\ with index $\alpha>0$} in the sense that
	\beam\label{eq:1}
	\P(X>x)= \dfrac{L(x)}{x^\alpha}\,,\qquad x>0\,,
	\eeam
	where $L$ denotes some \slvary\ \fct ; we refer to Bingham et al. 
	\cite{bingham:goldie:teugels:1987} for an encyclopedic treatment of
	univariate \regvary\ \fct s and to Resnick \cite{resnick:1987,resnick:2007}
	for the case of \regvary\ random vectors. 
	\par
	A natural question appears in this 
	context: given $Y$ is a non-negative \rv\ independent of 
	$X$, under which conditions is the product $X\,Y$ \regvary\ with index $\alpha$?
	This is a natural problem indeed: in numerous contexts of applied \pro y one 
	studies models which involve products of independent \rv s. Among those are 
	classical \ts\ models such as the {\em ARCH-GARCH family} and the {\em \sv\ model}; see Andersen et al. \cite{andersen:davis:kreiss:mikosch:2009} for an extensive treatment of
	these models in financial \tsa .
	In both cases, the real-valued time series $(X_t)$ is given via the relation
	$X_t=\sigma_t\,Z_t$, where $(\sigma_t)$ is a strictly stationary \seq\ of 
	positive \rv s which is either predictable \wrt\ the natural filtration of the iid 
	\seq\ $(Z_t)$ (such as for ARCH-GARCH) or $(\sigma_t)$ and $(Z_t)$ are
	mutually independent (such as for the \sv\ model). In both cases, there is strong interest in the 
	tail behavior of the products $X_t=\sigma_tZ_t$ (notice that, under the aforementioned conditions, $\sigma_t$ and $Z_t$  are independent). In the ARCH-GARCH
	the condition  $\E[|Z|^\alpha]<\infty$ and the dynamics of the volatlity \seq\ $(\sigma_t)$ ensure that
	$\P(\sigma_t>x)\sim c\,x^{-\alpha}$ for some positive constants $c,\alpha$ (for more details we refer the reader to Section~\ref{sec:RDE}).
	In turn, the condition $\E[|Z|^\alpha]<\infty$ and the so-called Breiman lemma 
	imply that 
	\beam\label{eq:holds}
	\P(\pm \sigma_tZ_t>x)\sim \E[(Z_t^{\pm})^\alpha]\,\P(\sigma_t>x)\,,\qquad\xto\,.
	\eeam 
	Breiman's result \cite{breiman:1965} is contained in the following
	useful lemma; for a proof, see Appendix C.3 in \cite{buraczewski:damek:mikosch:2016}.
	
	\ble\label{lem:breiman} Assume $X,Y$ are independent non-negative \rv s, $X$ is \regvary\ with 
	index $\alpha>0$ in the sense of \eqref{eq:1}, and $\E[Y^{\alpha+\delta}]< \infty$
	for some $\delta>0$ or $\P(X>x)\sim c\,x^{-\alpha}$ for some positive $c>0$
	and $\E[Y^\alpha]<\infty$.
	Then $\P(XY>x)\sim \E[Y^\alpha]\,\P(X>x)$ as $\xto$.
	\ele 
	
	Thus the \regvar\ of $X$ is preserved under multiplication
	with an independent non-negative \rv\ $Y$ if the corresponding assumptions
	on $Y$ hold, ensuring that $Y$ has a lighter tail than $X$. We already mentioned the case of an ARCH-GARCH process
	$(X_t)$ when $\sigma_t$ is \regvary\ with index $\alpha>0$
	and $X_t$ inherits this property if $\E[|Z_t|^{\alpha}]<\infty$. 
	In the \sv\ model, $X_t$ is \regvary\
	with index $\alpha>0$ if either $\sigma_t$ is \regvary\ with the same index
	and  $\E[|Z_t|^{\alpha+\delta}]<\infty$ for some $\delta>0$ and then 
	\eqref{eq:holds} holds,  or $Z_t$ is \regvary\ with index 
	$\alpha$, satisfying the {\em tail balance condition}\ :
	\beam\label{eq:tailbalance}
	\P(Z_t>x)\sim p_+ \dfrac{L(x)}{x^\alpha}\quad\mbox{and}\quad \P(Z_t<-x)\sim 
	p_-\dfrac{L(x)}{x^\alpha}
	\eeam
	for constants $p_\pm$ such that $p_++p_-=1$ and a \slvary\ \fct\ $L$, and
	$\E[\sigma_t^{\alpha+\delta}]<\infty$ for some $\delta>0$, and then
	\beao
	\P(\pm X_t>x)\sim \E[\sigma_t^\alpha] \,\P(\pm Z_t>x)\,,\qquad \xto\,,
	\eeao
	holds. 
	\par
	We mention that power-law tail behavior of a stationary \seq\ $(X_t)$
	is essential for the \asy\ behavior of their extremes and partial sums, and 
	related \pp\ \con\ and \fct als acting on them.
	For example, if $(Z_t)$ is iid and \regvary\ with index $\alpha>0$,
	then the \seq\ of the maxima $(a_n^{-1}M_n)$, where $M_n=\max_{i=1,\ldots,n} Z_i$,
	and $(a_n)$ satisfies $n\,\P(Z>a_n)\to 1$, converges in \ds\ to a
	Fr\'echet \ds\ $\Phi_\alpha(x)=\exp(-x^{-\alpha})$, $x>0$; see Embrechts et al.
	\cite{embrechts:klueppelberg:mikosch:1997}, Section~3.3. Moreover,
	the process of the points $(a_n^{-1}X_i)_{i=1,\ldots,n}$ converges in \ds\
	to an inhomogeneous Poisson process on $(0,\infty)$ with intensity \fct\
	$\alpha x^{-\alpha-1} dx$; see  
	 Resnick \cite{resnick:1987,resnick:2007}, Embrechts et al.
	\cite{embrechts:klueppelberg:mikosch:1997}, Chapter 5.
	Similarly, if $\alpha\in (0,2)$ and $Z$ is \regvary\ in the sense of 
	\eqref{eq:tailbalance} then $(a_n^{-1}(S_n-c_n))$ converges in \ds\
	(with suitable centering constants $(c_n)$)  
	to an infinite variance $\alpha$-stable limit; see Feller \cite{feller:1971}
	or Resnick \cite{resnick:2007}. Moreover, there is a vast literature that
	extends these results from the iid to the dependent case.
	\par
	For the completeness of presentation, we mention some related results for
	independent non-negative \rv s $X,Y$ when both are \regvary\ with the same
	index $\alpha$. This situation is much more subtle than the Breiman case.
	Still, $XY$ is \regvary\ with index $\alpha$: 
	
	\ble\label{lem:productregvaruniv}
		Assume that $X,Y$ are independent non-negative \rv s and  $X$ is \regvary\ with index $\alpha>0$. Then the following statements hold:
		\begin{enumerate}
			\item[\rm (1)] If either $Y$ is \regvary\ with index $\alpha$ or $\P(Y>x)=o\big(\P(X>x)\big)$ as $\xto$ 
				then $X\,Y$ is \regvary\ with index $\alpha$.
			\item[\rm (2)] If $\E[Y^\alpha]=\infty$ then 
				\beao
					\lim_{\xto}\dfrac{\P(X\,Y>x)}{\P(X>x)}= \infty\,.
				\eeao
			\item[\rm (3)] Assume that $X,Y$ are \regvary\ with index $\alpha>0$, $\E[X^\alpha+Y^\alpha]<\infty$,
				\beao
					c_0= \lim_{t\to\infty}\dfrac{\P(Y>t)}{\P(X>t)}\in [0,\infty)
				\eeao 
				and 
				\beam\label{eq:strange}
					\lefteqn{\lim_{M\to\infty}\limsup_{\xto}\dfrac{\P(X\,Y>x\,,M<X\le x/M)}{\P(X>x)}}\nonumber\\ 
					&=&\lim_{M\to\infty}\limsup_{\xto}\int_{M}^{x/M}\dfrac{\P(X>x/y)}{\P(X>x)}\,\P(Y\in dy)=0\,.
				\eeam
				Then 
				\beao
					\lim_{t\to\infty}\dfrac{\P(X\,Y>t)}{\P(X>t)}= \E[X^\alpha]+ c_0\,\E[Y^\alpha]\,.
				\eeao
		\end{enumerate}
	\ele
	The proof of this result is given in Appendix~\ref{app:prooflem1.2}.
	\bre 
	Condition \eqref{eq:strange} is a very technical assumption.
	To verify it one would need to have very precise information about the tail
	behavior of $X$. This condition does not follow from the uniform
	\con\ theorem for \regvary\ \fct s;
	the latter result ensures that for any $\vep>0$,
	\beao
	\lim_{\xto}\sup_{y\le \vep} \Big|\dfrac{\P(X>x/y)}{\P(X>x)}-y^\alpha\Big|=0\,.
	\eeao
	However, for the verification of  \eqref{eq:strange} we need 
	information about the deviation of $\P(X>x/y)/\P(X>x)$ from $y^\alpha$ 
	in the range  $y\in [M,x/M]$ for any $M>0$ and large $x$, 
	i.e., for large values of $y$. Part (3) was proved as Proposition~3.1
	by Davis and Resnick \cite{davis:resnick:1985a} in the case when $X,Y$
	are iid. In this case, \eqref{eq:strange} is necessary for
	$\P(XY>t)/\P(X>t)\to 2\E[X^\alpha]$ to hold.
	\ere

	We mention in passing that \regvar\ of $XY$ does in general not 
	imply \regvar\ of $X$ or $Y$; see Jacobsen et al.
	\cite{jacobsen:mikosch:rosinski:samorodnitsky:2009}.

\subsection{Closure of \regvar\ under multiplication -- the multivariate  case}
	Our main goal in this paper is to extend some of the aforementioned  results to the multivariate case. We start by introducing \regvar\ of random vectors. For this reason we equip $\R^{d_\bfX}$ with an arbitrary norm $\| \cdot \|$. A random vector $\bfX$ has a multivariate regularly varying distribution if $\| \bfX \|$ has a univariate regularly varying distribution 
and is asymptotically independent of  $\bfX / \| \bfX\|$ given $\|\bfX \|>x$. 
More precisely, we say that a random vector $\bfX\in \R^{d_\bfX}$ and its distribution are  {\em regularly varying} if
	\beam\label{eq:oo}
		\P \Big(\dfrac{\bfX}{\|\bfX\|}\in\cdot\,,\dfrac{\|\bfX\|}{x}\in \cdot\left| \|\bfX\|>x\Big)\right.  \stw\, \P (\bfTh_\bX\in\cdot)\,\P(Z\in\cdot)\,, \qquad \xto\,,
	\eeam
	where $Z$ is Pareto distributed with $\P(Z>y)=y^{-\alpha}$, $y>1$, and  $\bfTh_\bX$  assumes values in the unit sphere $\bbs^{d_\bfX-1}=\{\bfx\in\bbr^{d_\bfX}:\|\bfx\|=1\}$. The \ds\ of $\bfTh_\bX$ is the {\em spectral \ds } of $\bfX$.

	We will often refer to an equivalent formulation of multivariate regular variation. Namely, a random vector $\bfX\in \R^{d_\bfX}$ and its distribution are  {\em regularly varying} if and only if,
	 there exists a non-null Radon measure  $\mu^\bX$ on $\Ro^{d_\bfX} = \R^{d_\bfX} \setminus \{   \mathbf{0}\}$ such that
	\beao
		\mu^{\bX}_{t}(\cdot) = \dfrac{ \Prob ( t^{-1}\bfX \in  \cdot)}{\Prob(\| \bfX \| >t)} \vag \mu^\bX(\cdot)\,,\qquad t\to\infty\,,
	\eeao
	where $\stv$ denotes vague \con\ in the space of \ms s on $\Ro^{d_\bfX}$. Recall that for measures $\mu^{\bX}_{t}$, $\mu^\bX$ on $\Ro^{d_\bfX}$, $\mu^{\bX}_{t}  \vag \mu^\bX$ if for any function 
	$f$ from the set $C_c^+(\Ro^{d_\bfX})$ of  non-negative continuous functions on $\Ro^{d_\bfX}$ with compact support\footnote{ In the context of \regvar ,
the origin is excluded from consideration. Therefore
 a set $K\subset \Ro^{d_\bfX}$ is compact if it is compact in $\Ro^{d_\bfX}$ but bounded away from zero.} 
 we have
	\beao
		\int f(\bx) \: \mu_t^\bX(d \bx) \to \int f(\bx) \: \mu^\bX(d \bx)\,, \qquad t\to\infty\,.
	\eeao 
	It turns out that the limiting measure $\mu^\bX$ has the {\em homogeneity property.} 
	More precisely, there exists $\alpha_\bfX>0$ such that for any set $A$ in the Borel $\sigma$-algebra of $\Ro^{d_\bfX}$ we have
	\begin{equation*}
		\mu^\bX(t\,A) = t^{-\alpha_\bfX} \mu^\bX(A), \qquad t>0\,.
	\end{equation*} 
	We call $\alpha_\bfX$ the {\em index of regular variation} or {\em  tail index} of $\bfX$ and,  for short, we write $\bfX \in \RV(\alpha_\bfX,\mu^\bX)$.
	Of course, we necessarily have
	\beam\label{eq:L(x)}
		\Prob( \| \bfX \| >x) = \dfrac {L(x)}{x^{\alpha_\bfX}}\,,
	\eeam
	for some slowly varying function $L$.  We refer to Resnick \cite{resnick:1987,resnick:2007} as general references to multivariate \regvar\ and its applications.

\par
	
\par
Now consider two independent vectors 
$\bfX \in \RV(\alpha_\bfX,\mu^\bX)$ and $\bfY \in \RV(\alpha_\bfY,\mu^\bY)$ 
with values in $\bbr^{d_\bfX}$ and $\bbr^{d_\bfY}$, respectively.
Our goal is to establish sufficient conditions under which 
$\bfZ = \psi(\bfX, \bfY)$ is also regularly varying where
\beao
\psi \colon  \R^{d_\bfX} \times \R^{d_\bfY} \to \R^{d_\bfZ}
\eeao
is continuous, $a_\bfX$-homogeneous with respect to the 
first argument and $a_\bfY$-homogeneous with respect to the second one
for positive  $a_\bfX, a_\bfY$, i.e., for any $\bx \in \R^{d_\bfX}$ and 
	$\bfy \in \R^{d_\bfY}$,
	\beam\label{eq:homogen}
		\psi (s\,\bfx, t\,\bfy) = s^{a_\bfX}\,t^{a_\bfY}\, 
\psi(\bfx, \bfy)\,, \qquad  s,t\ge 0\,.	
	\eeam

	\bexam[\bf Products of independent regularly varying matrices]\rm
		If $d_\bfX = n_1 \cdot d_1$ 
then one can identify $\Ro^{d_\bfX}$ with the set of non-zero 
$n _1 \times d_1$ matrices $\mathbb{M}_{n_1 \times d_1}$. 
Similarly, if $ d_\bfY = d_1 \cdot m_1$, $\Ro^{d_\bfY} = \mathbb{M}_{d_1 \times m_1}$. 
	 We define 
$\psi ( \bfx , \bfy) = \bfx \cdot \bfy$ where $\bfx \cdot \bfy$ 
denotes ordinary matrix multiplication of an $n_1\times d_1$ 
matrix $\bfx$ with a $d_1\times m_1$ matrix $\bfy$. 
Then $d_\bfZ = n_1 \cdot m_1$, $a_\bfX=a_\bfY=1$, and 
$\bfZ$ is a product of two independent regularly varying matrices 
$\bfX$ and $\bfY$. 
\par
 In this case, regular variation of $\bfZ$ was proved in  
Basrak et al.~\cite{ basrak2002regular}; it is a  
		multivariate  analog of the Breiman 
Lemma~\ref{lem:breiman}: if
		\begin{equation*}
			 \bX \in \RV(\alpha,\mu^\bX) \quad \mbox{and} \quad   \E \big[\|\bY\|^{\alpha+\delta}\big]<\infty\quad\mbox{ for some $\delta>0$,}
		\end{equation*}
		then
		\beao
			\dfrac{\Prob(t^{-1}\bX \cdot \bY \in \cdot )}{\Prob(\| \bX \| >t )} \vag \   \eta(\cdot) := \Erw \big[\mu^\bX\big( \big\{\bx \: : \: \bx \cdot\bY \in \cdot\big\}\big)\big].
		\eeao
		In particular, if $\eta$ is non-null then $\bfZ = \bX \cdot \bY \in \RV \big( \alpha_\bfX,\mu^{\bfZ}\big)$ where
		\beao
			\mu^{\bfZ}(\cdot)=\dfrac{\eta(\cdot)} {\eta\big(\big\{ \bz\: : \: \|\bz\|>1 \big\}\big)}\,.
		\eeao
	\eexam

	\begin{example}[\bf Kronecker products of independent regularly varying matrices]\rm
		Suppose that $d_\bfX = n_1 \cdot d_1$ and 
$ d_\bfY = d_2 \cdot n_2$,  so we can identify $\Ro^{d_\bfX} = \mathbb{M}_{n_1 \times n_2}$, $\Ro^{d_\bfY} = \mathbb{M}_{d_1 \times d_2}$. Now define 
		$ \psi \colon \R^{d_\bfX} \times \R^{d_\bfY} \to \R^{n_1d_1n_2d_2}= 
\mathbb{M}_{n_1d_1 \times n_2d_2}$ via the Kronecker product $\psi(\bfx , \bfy) = \bfx \otimes \bfy$. As for ordinary matrix multiplication, we have $a_\bfX=a_\bfY=1$.
	\end{example}

	\begin{example}[\bf Random quadratic form]\rm 
		If $d_\bfY = d_\bfX^2$, identifying $\Ro^{d_\bfY} = \mathbb{M}_{d_\bfX \times d_\bfX}$, we define $\psi \colon \R^{d_\bfX} \times \R^{d_\bfY} \to \R$ 
by $\psi(\bfx, \bfy) = \bfx^\top \bfy \bfx$.  In this case, $a_\bfX=2$ and $a_\bfY=1$.
	\end{example}

\subsection{Organization of the article}
Our main result (Theorem~\ref{thm:main1}) yields sharp sufficient
conditions for \regvar\ of the homogeneous \fct\ $\psi(\bfX,\bfY)$
acting on independent \regvary\ random vectors $\bfX,\bfY$. 
The proof is given in Section~\ref{sec:proof1}.
We apply these results in Section~\ref{sec:RDE}. In particular,
in Section~\ref{subsec:prod}
we derive the \regvar\ properties of 
products of iid regularly varying quadratic matrices while,
in Section \ref{sec:RDE1}, we prove \regvar\ of solutions 
to affine \sre s under non-standard conditions.

\section{Main result}\setcounter{equation}{0}\label{sec:results}
 In what follows, $\bfX$ and $\bfY$ are independent \rv s with
values in $\bbr^{d_\bfX}$ and $\bbr^{d_\bfY}$, respectively, and we also
assume  $\bfX\in {\rm RV}(\alpha_\bfX,\mu^{\bfX})$ and 
 $\bfY\in {\rm RV}(\alpha_\bfY,\mu^{\bfY})$. We will study the \regvar\ property
of the $a_\bfX$-$a_\bfY$-homogeneous \fct\ $\bfZ=\psi(\bfX,\bfY)$; 
see \eqref{eq:homogen}.
We also need a {\em tail balance condition}\ : the following
limits exist and are finite  
	\begin{equation}\label{eq:balance}
		\lim_{t \to \infty}\dfrac{\P(\| \bX \|^{a_\bfX} > t)}{ 
\P (\| \bX \|^{a_\bfX} \cdot \|\bY\|^{a_{\bfY}} >t)} = c_\bX, \qquad 
\lim_{t \to \infty}\frac{\P(\| \bY \|^{a_\bfY} > t)}{ \P (\| \bX \|^{a_\bfX} \cdot \| \bY \|^{a_\bfY} >t)} = c_{\bY}\,.
	\end{equation}
We observe that $\|\bfX\|^{a_\bfX}$ and $\|\bfY\|^{a_\bfY}$ are 
\regvary\ with indices $\alpha_\bfX/a_\bfX$ and $\alpha_\bfY/a_\bfY$, respectively.
Therefore Lemmas~\ref{lem:breiman} and \ref{lem:productregvaruniv} apply:
\begin{itemize}
\item 
if $\alpha_\bfX/a_\bfX<\alpha_\bfY/a_\bfY$ then $\|\bfX\|^{a_\bfX}\|\bfY\|^{a_\bfY}$
is \regvary\ with index $\alpha_\bfX/a_\bfX$, $c_\bfX=1/\E[\|\bfY\|^{a_\bfY\,
\alpha_\bfX/a_\bfX}]\in 
(0,\infty)$ and $c_\bfY=0$.
\item
if  $\alpha_\bfX/a_\bfX=\alpha_\bfY/a_\bfY$ then  $\|\bfX\|^{a_\bfX}\|\bfY\|^{a_\bfY}$
is \regvary\ with index $\alpha_\bfX/a_\bfX$.
\item
if   $\alpha_\bfX/a_\bfX=\alpha_\bfY/a_\bfY$ and $\E [\|\bfY\|^{\alpha_\bfY}]=\infty$ 
then $c_\bfX=0$.
\item if $\alpha_\bfX/a_\bfX=\alpha_\bfY/a_\bfY$, 
$\E[\|\bfX\|^{\alpha_\bfX}+ \|\bfY\|^{\alpha_\bfY}]<\infty$, the limit 
\beam\label{eq:c_0}
\lim_{t\to\infty} \dfrac{\P(\|\bfY\|^{a_\bfY}>t)}{\P(\|\bfX\|^{a_\bfX}>t)}=c_0\in [0,\infty)
\eeam
exists and
\beam\label{eq:mm}
\lim_{M\to\infty}\limsup_{t\to\infty}\dfrac{\P(\|\bfX\|^{a_\bfX}\,\|\bfY\|^{a_\bfY}>t\,,M<{\|\bfX\|^{a_X}}\le t/M)}{\P(\|\bfX\|^{a_\bfX}>x)}=0
				\eeam
holds then 
\beam\label{eq:c_X}
c_\bfX= \dfrac{1}{\E[ \|\bY\|^{\alpha_\bfY}]+c_0\,\E[\|\bfX\|^{\alpha_\bfX}]}\qquad\mbox{ and }
c_\bfY=c_\bfX\,c_0\,.
\eeam
\end{itemize}

Now we formulate the first result of this paper.
\bth\label{thm:main1}
Assume that the $\bbr^{d_\bfX}$-valued 
$\bfX\in {\rm RV}(\alpha_\bfX,\mu^{\bfX})$ and the $\bbr^{d_\bfY}$-valued  
$\bfY\in {\rm RV}(\alpha_\bfY,\mu^{\bfY})$ random vectors are independent
and the balance condition \eqref{eq:balance} is satisfied for positive $a_\bfX,a_\bfY$.
Then the following relation holds for the $a_\bfX$-$a_\bfY$-homogeneous
\fct\ $\bfZ=\psi(\bfX,\bfY)$:
	\beam
	\lefteqn{\dfrac{\P \left(  t^{-1}\bfZ\in \cdot  \right)}{ \P (\| \bX \|^{a_\bfX} \cdot \| \bY \|^{a_\bfY} >t)}\stv \eta ( \cdot)}\nonumber\\
&= & \big(1 - c_\bX\,\E[ \| \bY\|^{(\alpha_\bfX a_\bfY)/a_\bfX}] - c_\bY\, 
\E[ \| \bX\|^{(\alpha_\bfY a_\bfX)/a_\bfY}]\big)
\times \E \big[\mu^\bX\big(\{ \bx\: : \: \psi(\bx, \Theta_\bY) \in \cdot \} \big)  \big]\nonumber \\  
	 	&& + c_\bX\E \left[ \mu^\bX(\{ \bx\: : \: \psi(\bx,  \bY) \in \cdot \})  \right] +  c_\bY \E \left[\mu^\bY(\{ \by\: : \:\psi(\bX, \by) \in \cdot \} ) \right].\label{eq:main}
	\eeam
	In particular, if $\eta$ is non-null, then $\bfZ  \in  \RV( \alpha_\bfZ, \mu^{\bfZ})$, where $\alpha_\bfZ = \frac{\alpha_\bfX}{a_\bfX}\wedge \frac{\alpha_\bfY}{a_\bfY}$ and
	\begin{equation*}
	\mu^{\bfZ}(\cdot) = \dfrac{\eta(\cdot)} {\eta\big(\big\{ \bz\: : \: \|\bz\|>1 \big\}\big)}\, .
	\end{equation*}
\ethe
Combining the discussion before Theorem~\ref{thm:main1} and the aforementioned results, we obtain the following con\seq es.

\bco\label{cor:cor1}
Assume the conditions of Theorem~\ref{thm:main1}.
	\begin{enumerate}
		\item[(1)] If $\frac{\alpha_\bfX}{a_\bfX} < \frac{\alpha_\bfY}{a_\bfY}$ then $c_\bfY=0$, $c_\bfX=1/\E[ \| \bY\|^{\alpha_\bfX a_\bfY/a_\bfX}]$, and 
\eqref{eq:main} holds with 
			\beao
				\eta(\cdot) = \dfrac{1}{\E[ \| \bY\|^{\alpha_\bfX a_\bfY/a_\bfX}] } \mu^\bX(\{\bx\: : \: \psi(\bx,\bY) \in \cdot \} )   \,.  
			\eeao
		\item[(2)] If $\P(\|\bX\|^{a_\bfX}>t) + 
\P(\| \bY \|^{a_\bfY}>t) = o (\P( \| \bX \|^{a_\bfX} \cdot \|\bY \|^{a_\bfY} > t))$ 
then $c_\bfX=c_\bfY=0$, and \eqref{eq:main} holds with
			\beao
				\eta(\cdot) = \E \left[ \mu^\bX( \{\bx\: : \: \psi(\bx, \Theta_\bY) \in \cdot \}) \right] \,.  
			\eeao
		
		\item[(3)] If $\alpha_\bfY/a_\bfY=\alpha_\bfX/a_\bfX$ and 
$\E\big[\|\bX\|^{\alpha_\bfX}+
\|\bY\|^{\alpha_\bfY }\big]<\infty$, \eqref{eq:mm} holds, and the limit $c_0$ in
\eqref{eq:c_0} exists, then $c_\bfX$ is given in \eqref{eq:c_X}, 
$c_\bfY=c_0\,c_\bfX$, and \eqref{eq:main} holds with
			\beam\label{eq:rr}
				\eta(\cdot)
				&=&\E\big[\mu^\bX\big(\{\bx: \psi(\bx, \bY) \in \cdot)\}\big)\big]+ c_0\,
				\E\big[\mu^\bY\big(\{\bfy: \psi(\bfX, \by) \in \cdot\}\big)\big]\,.
			\eeam
	\end{enumerate}
\eco
\bre As regards statement (2), one can verify that 
$\eta$ is symmetric  with respect to $\bX$ and $\bY$. In this case, necessarily
			$\frac{\alpha_\bfX}{a_\bfX}= \frac{\alpha_\bfY}{a_\bfY}$,  
and we can write
			\begin{align*}
				 \E \left[ \mu^\bX( \{\bx\: : \: \psi(\bx, \Theta_\bY) \in \cdot \}) \right] 
				 	= & \int_0^{\infty} \alpha_\bX r^{-\alpha_\bX-1} \P(\psi(r \Theta_\bX, \Theta_\bY) \in \cdot) \: dr \\
				 	= & \int_0^{\infty} \alpha_\bX r^{-\alpha_\bX-1} \P(\psi(\Theta_\bX, r^{a_\bX/ a_\bY}\Theta_\bY) \in \cdot) \: dr \\
				 	= & \int_0^{\infty} \alpha_\bY r^{-\alpha_\bY-1} \P(\psi(\Theta_\bX, r\Theta_\bY) \in \cdot) \: dr \\
				 	= &  \E \left[ \mu^\bY( \{\by\: : \: \psi(\Theta_\bX, \by) \in \cdot \}) \right].
			\end{align*} 
\ere
\section{Proof of Theorem~\ref{thm:main1}}\label{sec:proof1}\setcounter{equation}{0}	
	Throughout this section we consider an $\bbr^{d_\bfX}$-valued
$\bX \in \RV(\alpha_\bfX, \mu^\bfX)$ random vector independent of 
an $\bbr^{d_Y}$-valued $\bY \in \RV(\alpha_\bfY, \mu^\bfY)$. Recall that $\bfZ=\psi(\bfX,\bfY)\in\bbr^{d_\bfZ}$.
	Take any function $f$ from the set $C^+_c(\Ro^{d_\bfZ})$ of nonnegative continuous functions with compact support in $\Ro^{d_\bfZ}$.
Write 
\beao
\eta_t(\cdot) = \dfrac{\P(t^{-1}\bfZ\in\cdot)}{
\P(\|\bfX\|^{a_\bfX}\|\bfY\|^{a_\bfY}>t)}\,.
\eeao
Then  \eqref{eq:main} turns into $\eta_t\stv \eta$ as $t\to\infty$
which can be re-formulated as 
\beao	
\lim_{t \to \infty}\dfrac{\E[f( t^{-1} \bfZ)]}{ \P (\| \bX \|^{a_\bfX} \cdot 
\|\bY\|^{a_\bfY} >t)}= \lim_{t\to\infty} \int f(\bfz) \eta_t(d\bfz) 
 = \int f(\bz) \: \eta(d\bz)\,,\quad f\in  C^+_c(\Ro^{d_\bfZ})\,.
	\eeao
Since $\psi$ is continuous
	\begin{equation*}
		M_\psi = \sup \{ \| \psi(\bfx, \bfy) \| \: : \:  \|\bfx\|=1, \: \| \bfy\| =1 \}<\infty.
	\end{equation*}
It is also $a_\bfX$-$a_\bfY$-homogeneous and therefore
	\begin{equation*}
		\| \psi(\bfx, \bfy)\| \leq M_\psi \|\bfx\|^{a_\bfX} \|\bfy\|^{a_\bfY}. 
	\end{equation*}
Then we also have for any set $A_r=\{\bfz: \|\bfz\|>r\}$, $r>0$,
in view of \regvar\ of $\|\bfX\|^{a_\bfX}\|\bfY\|^{a_\bfY}$,
\beao
\sup_{t>0}\eta_t( A_r)\le \dfrac{\P\big( M_\psi \|\bfX\|^{a_\bfX} \|\bfY\|^{a_\bfY}>r\,t \big))}{\P(\|\bfX\|^{a_\bfX}\|\bfY\|^{a_\bfY}>t)}<\infty\,.
\eeao
It follows from Resnick \cite{resnick:1987}, Proposition~3.16,
that $(\eta_t)$ is vaguely relatively compact. Hence $(\eta_{t_k})$
converges vaguely along \seq s $t_k\to\infty$ as $\kto$, and it remains
to show that these limits coincide with $\eta$.
\par
The proof of the theorem is given through several auxiliary result
which we provide first. The main steps of the proof are given at
the end of this section.\\[2mm]
{\bf Limits of $\E[f \big( t^{-1} \psi(\bX,\bY)\big)\mid\bfY]$.} 
By \regvar\ of $\bfX$ we have
\beam\label{eq:regvarX}
\mu_{t}^\bfX( \cdot)=  \dfrac{\P(t^{-1} \bfX\in\cdot)}{\P(\|\bfX\|>t)}\stv \mu^\bfX(\cdot)\,,\qquad t\to\infty\,.
\eeam
Define
\begin{equation}\label{eq:defGt}
		g_t(\by) = 
\dfrac{\E \left[ f\left(t^{-1}\psi(\bX,\by)\right)\right]}{\P ( \| \bX \|^{a_\bfX} >t )} = \int f(\psi( \bx, \by))  \: \mu_{t^{1/a_\bfX}}^\bX(d\bx) , \qquad \by \in 
\R^{d_\bfY}, \: t>0.
	\end{equation}
In view of  \eqref{eq:regvarX} we expect that the \rhs\ converges to
	\beam\label{eq:aa}
g_t(\bfy)\to		g(\by) = \int f(\psi(\bx, \by))  \: \mu^\bX(d\bx)
<\infty\,, \qquad t\to\infty \,, \qquad \by \in \R^{d_\bfY} . 
	\eeam
However, the \fct\ $\bx \mapsto f(\psi(\bx,\by))$ may not have 
compact support and therefore some additional argument is needed.
\ble \label{lem:convergence}
Relation \eqref{eq:aa} holds for any $f \in C^+_c(\Ro^{d_\bfZ})$. 
\ele
\begin{proof}[Proof of Lemma~\ref{lem:convergence}]
Fix $\by \in \R^{d_\bfX}$. Since $f$ is compactly supported there are constants
 $M_f, c_f>0$ such that
\begin{equation}\label{eq:suppf}
			 {\rm supp} (f) \subseteq \{ \bz \in \R^m \: : \: c_f^{-1} \leq \| \bz \| \leq c_f \}\quad\mbox{and}\quad \sup_{\bz \in \R^{d_\bfZ}}f(\bz) \leq M_f.
		\end{equation} 
For $r\geq 1$ choose any continuous \fct\ 
$\varphi_r \colon \R^{d_\bfX} \to [0,1]$ such that 
\beao
\varphi_r(\bx) = \left\{\barr{ll}
1\,,& \|\bx \| \leq r\,,\\
0\,,&  \| \bx \| \geq 2r\,. \earr\right.
\eeao
We have
		\begin{equation*}
			g_t(\by) = \int f(\psi(\bfx,\bfy))\,\varphi_r(\bx)  \: \mu^\bX_{t^{1/{a_\bfX}}}(d\bx) +\int f(\psi(\bx, \by))(1-\varphi_r(\bx))  \: 
\mu^\bX_{t^{1/{a_\bfX}}}(d\bx)=I_1+I_2\,.
		\end{equation*}
		The contribution of the second term is negligible 
since in view of \eqref{eq:regvarX},
\beao
0&\le& \lim_{r\to\infty}\limsup_{t\to\infty} I_2\\
&\le & M_f \,\lim_{r\to\infty}\lim_{t\to\infty}
\mu^\bX_{t^{1/a}}( \{  \bx \: : \: \|\bx\| >r \})\\
&=& M_f \,\lim_{r\to\infty}
\mu^\bX( \{  \bx \: : \: \|\bx\| >r \})=0\,.
\eeao
Thus it suffices to prove $\lim_{r \to \infty} \lim_{t \to \infty} I_1 = g(\by)$.
The function $\bx \mapsto  f(\psi(\bx, \by))\varphi_r(\bx) $ is 
continuous and non-negative for any choice of $\by \in \R^{d_\bfY}$ and $r>1$, and 
its support is contained in $\{ \bx \in \R^{d_\bfX} \: : \: (M_\psi\| \by\|^{a_\bfY} c_f)^{-1/a_{\bfX}} \leq \| \bx \| \leq 2r\}$ which is a compact subset of $\Ro^{d_\bfX}$. Regular variation of $\bX$ and monotone \con\ allow one to take the successive limits
\beao
 \lim_{r\to\infty}\lim_{t\to\infty} I_1&=&\lim_{r\to\infty}\int f(\psi(\bx, \by))\varphi_r(\bx)  \: \mu^\bX(d\bx)= g(\by)\\
&=&
\int_{\|\bx\| \geq (M_\psi\| \by\|^{a_\bfY} c)^{-1/a_\bfX} } f(\psi(\bx, \by)) \: \mu^\bX(d\bx) \leq M_f (M_\psi\| \by\|^a_\bfY c)^{\alpha_\bfX/a_\bfX}<\infty\,.
\eeao
	\end{proof}	
 	The next result presents a  continuity bound for $g_t$.
	
	  \ble \label{lem:2}
	  	Let $ f \in C^+_c(\R^{d_\bfZ}_{ \bf0})$. 
For any $\varepsilon>0$ one can choose  $\delta>0$ and $t_0>0$
such that for any $\bs$, $\br \in \bbs^{d_\bfX-1}$ with $\|\bs-\br\| \leq \delta$
and any $t>t_0$,
	  	\beam\label{eq:sun1}
	  		|g_t(\br) - g_t(\bs)| \leq \varepsilon\,.
	  	\eeam
	  \ele
	 
	  \begin{proof}[\bf Proof of Lemma~\ref{lem:2}]
	  	Fix $\varepsilon_1>0$. Choose $M_f, c_f>0$ from \eqref{eq:suppf}. By uniform continuity of $f$ we can choose  $\eta \in (0, \varepsilon_1)$ such that
	  	$\|\bz_1 - \bz_2\| \leq \eta$ implies $\| f(\bz_1 )- f(\bz_2 ) \| \leq \varepsilon_1$. 
	  	Since $\psi$ is uniformly continuous on $\bbs^{d_\bfX-1}\times 
\bbs^{d_\bfY-1}$ we can find $\delta >0$ such that for $\br, \bs \in \bbs^{ d_\bfY-1}$ with  $\| \br - \bs \| < \delta$,
	  	\begin{equation*}
	  		 \| \psi(\bfx, \br) -\psi(\bfx, \bs)\| < \eta^{2}, \qquad \|\bfx\| =1.
	  	\end{equation*}
	  	Then by homogeneity of $\psi$,
	  	\begin{equation*}
	  		 \| \psi(\bfx, \br) -\psi(\bfx, \bs)\| < \| \bfx\|^{a_\bfX} \eta^{2}, \qquad \bfx\in \R^{d_\bfX}\,,
	  	\end{equation*}
                and we can write for $t>0$, 
	  	\begin{align*}
	  		|g_t(\bs) - g_t(\br) | \leq &\int |f(\psi(\bx,\bs)) - f(\psi(\bx,\br)) | \: \mu^\bX_{t^{1/a_\bfX}}(d \bx)\\
	  					 =& \int_{\|\bx\| \geq (M_\psi c_f)^{-1/a_\bfX}}  |f(\psi(\bx,\bs)) - f(\psi(\bx,\br)) | \: \mu^\bX_{t^{1/a_\bfX}}(d \bx) \\
	  					 \leq &\int_{\| \bx  \|>\varepsilon_1^{-1/a_\bfX} }| f(\psi(\bx,\bs)) - f(\psi(\bx,\br)) |\: \mu^\bX_{t^{1/a_\bfX}}(d \bx) \\
	  						&+  \int_{\| \bx  \| \leq \eta^{-1/a_\bfX}, \:\|\bx\| \geq (M_\psi c_f)^{-1/a_\bfX}}  |f(\psi(\bx,\bs)) - f(\psi(\bx,\br))| \: \mu^\bX_{ t^{1/a_\bfX}}(d \bx) \\
						 \leq& 2\, M_f\, \varepsilon_1^{\alpha_\bfX/a_\bfX}   \frac{L(\varepsilon_1^{-1/a_\bfX}t^{1/a_\bfX})}{L(t^{1/a_\bfX})} + 
\varepsilon_1 (M_\psi c)^{\alpha_\bfX/a_\bfX},
	  	\end{align*}
	  	where $L$ is defined in \eqref{eq:L(x)}.
	  	Given $\varepsilon>0$ we can choose $\varepsilon_1$ sufficiently small such that
	  	\begin{equation*}
	  		2\, M_f\, \varepsilon_1^{\alpha_\bfX/a_\bfX}   \frac{L(\varepsilon_1^{-1/a_\bfX}t^{1/a_\bfX})}{L(t^{1/a_\bfX})} + 
\varepsilon_1 (M_\psi c)^{\alpha_\bfX/a_\bfX} \leq \frac \varepsilon 3  \frac{L(\varepsilon_1^{-1/a_\bfX}t^{1/a_\bfX})}{L(t^{1/a_\bfX})}  +  \frac \varepsilon 3.
	  	\end{equation*}
	  	Choosing $t_0$ big enough,  one ensures that
	  	\begin{equation*}
	  		 \frac{L(\varepsilon_1^{-1/a_\bfX}t^{1/a_\bfX})}{L(t^{1/a_\bfX})} \leq 2
	  	\end{equation*} 
	  	which proves the claim.
	  \end{proof}

Note that by continuity of $f$ and $\psi$,  $g$ is also continuous on $\R^{d_\bfX}$, hence also uniformly continuous  on the unit sphere $\bbs^{d-1}$.
We will use this comment in the proof of the next lemma.

	\ble \label{cor:unif}
		Let $ f \in C^+_c( \R^{d_\bfZ}_{\bf0})$. Then $g_t\to g$ as $t\to\infty$
 uniformly on $\bbs^{ d_\bfY-1}$.
	\ele
	\begin{proof}[\bf Proof of Lemma~\ref{cor:unif}.]
		Fix $\varepsilon > 0$ and take $\delta>0$, $t_0>0$ that  satisfy the claim of Lemma~\ref{lem:2} and
		\begin{equation*}
			\| \bs - \br \| \leq \delta \Rightarrow | g(\br) - g(\bs)| \leq \varepsilon.
		\end{equation*}
Let $\{ \br_k\}_{k=1}^N$ for $N = N(\delta)$ be a
		$\delta$-covering of $\bbs^{ d_\bfY-1}$. Take $t_1>0$ so large that
		\begin{equation*}
		 	\max_{1\leq k \leq N}| g_t(\br_k) - g(\br_k)| \leq \varepsilon, \quad t>t_1.
		\end{equation*}
		Then for any $\bs \in \bbs^{d_\bfY-1}$ we have $\|\bs - \br_k\|\leq \delta$ for some $k$ and for $t> t_0\vee t_1$ we have
		\begin{equation*}
			|g_t(\bs) - g(\bs)| \leq |g_t(\bs) - g_t(\br_k)| +  |g_t(\br_k) - g(\br_k)| + |g(\br_k) - g(\bs)| \leq 3\varepsilon.
		\end{equation*}
This finishes the proof of the lemma.
	\end{proof}

	Before we proceed with the final steps in the proof of 
Theorem~\ref{thm:main1} we observe that homogeneity of $\mu^\bX$ and $\psi$ 
implies for any $r>0$ and $\by \in \R^{d_\bfY}$,
	\begin{equation*}
		g(r\,\by) = r^{\frac{\alpha_\bfX a_\bfY}{a_\bfX}}\,g(\by).
	\end{equation*}
	Now we define functions $h_t \colon \R^{d_\bfX} \to [0, +\infty)$ by 
	\begin{equation*}
		h_t(\bx) = \int f(\psi(\bx, \by))  \: \mu^\bY_{t^{1/a_\bfY}}(d\by) = \frac{\E \left[ f\left(t^{-1}\psi(\bx,\bY)\right)\right]}{\P ( \| \bY \|^{a_\bfY} >t )}, \qquad \bx \in \R^{d_\bfX}, \: t>0.
	\end{equation*}
By a symmetry argument, interchanging the roles of $\bY$ and $\bX$,
	we conclude that $h_t\to h$ as $t\to\infty$ point-wise 
in $\bbr^{d_\bfX}$ and   uniformly  on $\bbs^{d_\bfX-1}$ where 
	\begin{equation*}
		h(\bx) = \int f(\psi(\bx, \by))  \: \mu^\bY(d\by), \qquad \bx \in \R^{d_\bfX}\, .
	\end{equation*}
	The limiting function  is also homogeneous, i.e., for $r>0$ and $\bx \in \R^{d_\bfX}$,
	\begin{equation*}
		h(r\,\bx) = r^{\frac{\alpha_\bfY a_\bfX}{a_\bfY}}\,h(\bx).
	\end{equation*}

\begin{proof}[\bf Main steps in the proof of Theorem~\ref{thm:main1}]
	Recalling the notation introduced so far, our goal is to 
prove \eqref{eq:main} in disguised form by applying an approach via test \fct s:
	\beao
	\lefteqn{	\lim_{t \to \infty}\dfrac{\E \left[ f \left( t^{-1}\psi(\bX, \bY) \right) \right]}{ \P (\| \bX \|^{a_\bfX} \cdot \| \bY \|^{a_\bfY} >t)}}\\&= &
\big(1 - c_\bX\,\E\big [\| \bY\|^{\alpha_\bfX a_\bfY/a_\bfX}\big]- c_\bY\, 
\E\big[ \| \bX\|^{a_\bfX\alpha_\bfY/ a_\bfY}\big)\, 
		\E \left[g( \Theta_\bY)\right] \\&  & + c_\bX\,\E \left[g(\bY)\right] +  c_\bY \,\E \left[h(\bX) \right]\,.
\eeao
	Choose $M_f>0$ from \eqref{eq:suppf} and  consider the following decomposition, for $\eta \in (0,1)$, 
	\begin{align*}
		\E \left[ f \left( t^{-1} \psi(\bX, \bY) \right) \right]=  
&  \E \left[ f \left( t^{-1}\psi(\bX,\bY) \right) \1\big( \| \bY \|^{a_\bfY} \leq \eta t\big) \right] \\
									& + \E \left[ f \left( t^{-1}\psi(\bX,\bY) \right) \1\big(\| \bX \|^{a_\bfX} \leq \eta t, \: \| \bY \|^{a_\bfY} > \eta t\big) \right] \\
									& + \E \left[ f \left( t^{-1}\psi(\bX, \bY) \right) \1\big(\| \bX \|^{a_\bfX} > \eta t, \: \| \bY \|^{a_\bfY} > \eta t\big) \right] \\
									= & J_1(t) + J_2(t) + J_3(t). 
	\end{align*}
	Since $f$ is bounded and $\bfX,\bfY$ are independent we have 
$J_3(t) = o(\P (\| \bX \|^{a_\bfX} \cdot \| \bY \|^{a_\bfY} >t))$. Thus it remains
 to investigate $J_1$ and $J_2$. We begin with the analysis of the first term, since it requires more work.\\
{ \bf Analysis of $J_1$.}	We claim that
	\beao
\lefteqn{	\lim_{\eta \to 0}\liminf_{t \to \infty}\frac{J_1(t)}{ \P (\| \bX \|^{a_\bfX} \cdot \|\bY\|^{a_\bfY} >t)}}\\& = & \lim_{\eta \to 0}\limsup_{t \to \infty}  \frac{J_1(t)}{ \P (\| \bX \|^{a_\bfX} \cdot \|\bY\|^{a_\bfY} >t)} \\ &= &  
(1 - c_\bX\,\E \| \bY\|^{\alpha_\bfX a_\bfY/ a_\bfX} - c_\bY\, \E \| \bX\|^{\alpha_\bfY a_\bfX/ a_\bfY}) \E \left[ g(\Theta_\bY)\right]+  c_\bX\E \left[ g(\bY)\right].
	\eeao
	Below we will present a detailed argument for
	\beam\label{eq:Jsup}
	\lefteqn{\lim_{\eta \to 0}\limsup_{t \to \infty}  \frac{J_1(t)}{ \P (\| \bX \|^{a_\bfX} \cdot \|\bY\|^{a_\bfY} >t)}}\nonumber\\& \leq& (1 - c_\bX\,\E \| \bY\|^{\alpha_\bfX a_\bfY/ a_\bfX} - c_\bY \E \| \bX\|^{\alpha_\bfY a_\bfX/ a_\bfY}) \E \left[ g(\Theta_\bY)\right]+  c_\bX\E \left[ g(\bY)\right].\nonumber\\
	\eeam
	The lower bound can be established in a similar fashion. Write for $\bfz\ne \bf0$, $\wt \bfz=\bfz/\|\bfz\|$, and 
	\begin{align*}
		J_1(t) & = \int_{\| \by \|^{a_\bfY} \leq \eta t} \E \left[ f \left( t^{-1}\psi(\bX,\by) \right) \right] \: \P (\bY \in d \by)\\
			& = \int_{\| \by \|^{a_\bfY} \leq \eta t}  g_{\frac{t^{1/a_\bfX}}{\| \by \|^{a_\bfY/a_\bfX}}} \left(\wt \by \right) \P( \| \bX \|^{a_\bfX} \cdot \|\by \|^{a_\bfY} > t ) \: \P (\bY \in d \by),
	\end{align*}
	where $g_t$ is given via~\eqref{eq:defGt}. By virtue 
of Lemma~\ref{cor:unif}, for any $\varepsilon>0$ there is a sufficiently 
small $\eta>0$ such that
	\beao
\lefteqn{\left| J_1(t) - \int_{\| \by \|^{a_\bfY} \leq \eta t}  
g \left(\wt \bfy  \right) \P( \| \bX \|^{a_\bfX} \cdot \|\by \|^{a_\bfY} > t ) \: \P (\bY \in d \by) \right|}\\
&\le &\int_{\| \by \|^{a_\bfY} \leq \eta t} \Big|g_{\frac{t^{1/a_\bfX}}{\| \by \|^{a_\bfY/a_\bfX}}} \left(\wt \by \right) -g \left(\wt \by \right)\Big|
\;\P( \| \bX \|^{a_\bfX} \cdot \|\by \|^{a_\bfY} > t ) \: \P (\bY \in d \by)\\[2mm]& \leq& \varepsilon \, \P( \|\bX\|^{a_\bfX} \cdot \|\bY\|^{a_\bfY} > t ).
	\eeao
	Thus, since $\varepsilon$ is arbitrary, 
we only need to investigate the expectation
	\beao
I(t)=\E \left[ g (  \wt\bY) 
\;\1\big(\|\bX\|^{a_\bfX}\cdot \|\bY \|^{a_\bfY} > t\,, \|\bY \|^{a_\bfY} \leq \eta t\big) \right].
\eeao

		If $\E[g(\wt \bfY)]=0$ then  by homogeneity of $g$, 
$g(\bY) =0 $ a.s. which implies $\E [g (\bY)] = 0$
		and $\E[g(\Theta_\bY)] =0$,
		so the claim follows trivially. Now assume $\E[g(\wt \bfY)]>0$.  
		Let $Y'$ be a~random variable independent of 
$\bX$ and $\bY$ with distribution given by
		\begin{equation*}
			\P(Y' \in \cdot ) =  
\E\Big[ \dfrac{g (  \wt \bY)}{\E[g(\wt \bfY)]}\1\big(\| \bY \|^{a_\bfY} \in \cdot\big)\Big].
		\end{equation*}
		Then, by regular variation of $\bY$, as $t\to\infty$,
		\begin{equation*}
			\dfrac{\P(Y'>t)}{\P(\| \bY \|^{a_\bfY}>t)} = 
\E \Big[   \dfrac{g (  \wt \bY)}{\E[g(\wt \bfY)]}\,\Big|\,\|\bY \|^{a_\bfY} > t \Big] \to \dfrac{\E[g(\Theta_\bY)]}{
\E[g(\wt \bfY)]}\,.
		\end{equation*}
		Therefore for any $\delta>0$ there exists $T=T(\delta)$ 
such that 
		\beam\label{eq:needed}
			\Big|\dfrac{\P(Y'>t)}{\P(\| \bY \|^{a_\bfY}>t)} -
\dfrac{\E[g(\Theta_\bY)]} 
{\E[g(\wt \bfY)]}\Big|\le\delta\,,\qquad t\ge T\,. 
		\eeam
		Without loss of generality 
we may assume that $T\uparrow \infty$ when $\delta\downarrow 0$. 
Consider the following decomposition
\beao
\dfrac{I(t)}{\E[g(\wt \bfY)]}&=&
\P( \| \bX \|^{a_\bfX}\, Y'>t, \: Y' \leq \eta t) \\
				&=&  \P( \| \bX \|^{a_\bX} Y'>t, \: Y' >T) 
				+  \P( \| \bX \|^{a_\bX} Y'>t, \: Y' \leq T)  
				-   \P( \| \bX \|^{a_\bX} Y'>t, \:  Y'>\eta t) \\
				&=&  I_1(t)+I_2(t)-I_3(t).
\eeao
By Breiman's Lemma~\ref{lem:breiman} and definition of $c_\bfX$ we have
\beao
\lim_{\delta\downarrow 0}\lim_{t\to\infty}\dfrac{\E[g(\wt Y)]\,I_2(t)}{\P(\|\bfX\|^{a_\bfX}\,\|\bfY\|^{a_\bfY}>t)}&=&\lim_{\delta\downarrow 0}\lim_{t\to\infty}\dfrac{\E[g(\wt Y)]\,I_2(t)}{\P\big(\|\bfX\|^{a_\bfX}>t\big)}
\dfrac{\P\big(\|\bfX\|^{a_\bfX}>t\big)}{\P(\|\bfX\|^{a_\bfX}\,\|\bfY\|^{a_\bfY}>t)}\\
&=&\lim_{\delta\downarrow 0}c_\bfX\,\E[g(\wt Y)]\;
\E[(Y')^{a_\bfX/\alpha_X}\,\1(Y'\le T)]\\
&=&c_\bfX\,\E \big[g (  \wt \bY)  \|\bY \|^{\alpha_\bfX a_\bfY/a_\bfX} \big]\\
&=&c_\bfX\,\E \big[g (\bY)\big] \,.
\eeao
For the first term we have by \eqref{eq:needed}\,,
		\beao
		\E[g(\wt Y)]\,	I_1(t)&=&\E[g(\wt Y)]\,
\int_T^\infty \P\big( Y'> T \vee (t/\|\bfx\|^{a_\bfX})\big) \,\P(\bfX\in d\bfx))\\
&\le & (1+\delta)\,\E[g(\Theta_\bfY)]\,\int_T^\infty \P\big( \|\bfY\|^{a_\bfY}> T \vee (t/\|\bfx\|^{a_\bfX})\big) \,\P(\bfX\in d\bfx))\\
&=& (1+\delta)\E[g(\Theta_\bfY)]\,\P\big( \| \bX \|^{a_\bfX} \| \bY \|^{a_\bfY} >t, \: \| \bY \|^{a_\bfY}>T \big)\\
&=& (1+\delta)\,\E[g(\Theta_\bfY)]\,\Big[ \P( \| \bX \|^{a_\bfX} \| \bY \|^{a_{\bfY}} >t )\\&& - \P( \| \bX \|^{a_\bfX} \| \bY \|^{a_\bfY} >t, \: \| \bY \|^{a_\bfY} \leq T ) \Big] \\
&\sim & (1+\delta)\,\E[g(\Theta_\bfY)]\,
\P\big( \| \bX \|^{a_\bfX }\, \| \bY \|^{a_\bfY} >t \big)\\&&\times \Big[1  -
 \E\big[\|\bY \|^{\alpha_\bfX a_\bfY/a_\bfX} \1\big(\| \bY \|^a_\bfX \leq T\big)\big]\,
\dfrac{ \P(\| \bX \|^{a_\bfX} > t)}{\P\big( \| \bX \|^{a_\bfX }\, \| \bY \|^{a_\bfY} >t \big)}\Big] \,.
		\eeao
In the last step we used Breiman's result as $t\to\infty$. 
Now, recalling the definition of $c_\bfY$, we conclude that
\beao \lefteqn{
\lim_{T\to\infty}\limsup_{t\to\infty}
\dfrac{	\E[g(\wt Y)]\,	I_1(t)}{\P\big( \| \bX \|^{a_\bfX }\, \| \bY \|^{a_\bfY} >t \big)}}\\
&\le &  (1+\delta)\,\E[g(\Theta_\bfY)]\,
 \Big[1  - c_\bfX\,
 \E\big[\|\bY \|^{\alpha_\bfX a_\bfY/a_\bfX}\big]\Big]\,,
\eeao
and the corresponding lower bound can be derived in an analogous way for any small
$\delta>0$.
\par
Finally, we deal with the third term. 
First we observe that, by \regvar ,
		\beam \label{eq:miniclaim}
\lefteqn{\lim_{\eta \downarrow 0} 
\lim_{t \to \infty}\dfrac {\P\left (\| \bX \|^{a_\bfX} > \eta^{-1}\,, \:  
\| \bY \|^{a_\bfY} > \eta\, t \right) }{\P\left (\| \bX \|^{a_\bfX} \cdot 
\| \bY \|^{a_\bfY} > t \right)}}\nonumber\\&=&\lim_{\eta \downarrow 0} 
\P\big(\| \bX \|^{a_\bfX} > \eta^{-1}\big)\,
\eta^{-\alpha_\bfY/a_\bfY}\,
\lim_{t \to \infty}\dfrac {\P\big ( \| \bY \|^{a_\bfY} >  t \big) }{\P\left (\| \bX \|^{a_\bfX} \cdot 
\| \bY \|^{a_\bfY} > t \right)}\nonumber\\
&=& c_\bfY\,\lim_{\eta \downarrow 0} 
\P\big(\| \bX \|^{a_\bfX} > \eta^{-1}\big)\,
\eta^{-\alpha_\bfY/a_\bfY}\,
=0.
		\eeam		
		Indeed, if $\E[ \|\bX\|^{\alpha_Y a_\bfX/a_\bfY}] = \infty$ 
then $c_\bY=0$ and therefore the \rhs\ is zero; see Lemma~\ref{lem:productregvaruniv}(2).
		On the other hand, if  
$\E[ \|\bX\|^{\alpha_\bfY a_\bfX/a_\bfY}] < \infty$ then 
\beao
\P\big(\| \bX \|^{a_\bfX} > \eta^{-1}\big)&=&
\P\big( \| \bX \|^{a_\bfX\alpha_Y/a_\bfY} > \eta^{-\alpha_Y/a_\bfY}\big)
=o( \eta^{\alpha_\bfY/a_\bfY})\,,\qquad
\eta\downarrow 0\,, 
\eeao 
and therefore the \rhs\ in \eqref{eq:miniclaim} is zero.
\par
With \eqref{eq:needed} and Breiman's result  
at hand, we have as $t\to\infty$,
\beao
\lefteqn{\E[g(\wt Y)]	\,I_3(t)}\\ & \leq & (1+\delta) \E[g(\Theta_\bfY)]\,
\P\big( \| \bX \|^{a_\bfX}\, \| \bY \|^{a_{\bfY}} >t , \| \bY \|^{a_\bfY} > \eta\, t)\\
		&	  =&  (1+\delta)\,\E[g(\Theta_\bfY)]\,\\&&\times
\Big[\P\left ( \| \bX \|^{a_\bfX} > \eta^{-1} , \| \bY \|^{a_\bfY} > \eta\, t \right)
+  \P\left ( \| \bX \|^{a_\bfX} \leq \eta^{-1} , \| \bX \|^{a_\bfY} \cdot 
\| \bY \|^{a_\bfY} > t \right)\Big] \\ &\sim &  (1+\delta)\,\E[g(\Theta_\bfY)]\,\Big[
\P\big ( \| \bX \|^{a_\bfX} > \eta^{-1}\, , 
\| \bY \|^{a_\bfY} > \eta\, t \big)\\&& +  \E\big[ \|\bX\|^{a_\bfX \alpha_\bfY/a_\bfY}\,
 \1\big(\| \bX \|^{a_\bfX} \leq \eta^{-1}\big)\big]\,
\P\big(\|\bfY\|^{a_\bfY}>t\big)\,,\qquad t\to\infty\,. 
\eeao
Now an application of \eqref{eq:miniclaim} and the definition of 
$c_\bfY$ yield
\beao
\lim_{\eta \downarrow 0}\limsup_{t\to\infty}
\dfrac{\E[g(\wt Y)]\,I_3(t)}{\P\big (\| \bX \|^{a_\bfX} \cdot 
\| \bY \|^{a_\bfY} > t \big)}\le c_\bfY\,(1+\delta)\,  \E[g(\Theta_\bfY)]\,
 \E\big[ \|\bX\|^{a_\bfX \alpha_\bfY/a_\bfY}\big]
\,.
\eeao
This establishes an upper bound; the corresponding lower bound is completely
analogous. This proves \eqref{eq:Jsup}.\\[2mm]
{\bf Analysis of $J_2$. } This term is significantly simpler since we have
\begin{equation*}
	J_2(t)= \int_{\| \bx \|^{a_\bfx} \leq \eta t}  h_{\frac{t^{1/a_\bfY}}{\| \bx \|^{a_\bfX/a_\bfY}}} \left(\wt \bx \right) \P \left( \| \bY \|^{a_\bfY} \cdot \left(\frac 1\eta \wedge \|\bx \|^{a_\bfX} \right) > t \right) \: \P (\bX \in d \bx).
\end{equation*} 
Appealing to dominated convergence theorem, we obtain
\begin{equation*}
	\lim_{t \to \infty} \frac{J_2(t)}{\P (\| \bfX \|^{a_{\bfX}}\|\bfY\|^{a_\bfY} >t   )} =  \int_{\R^{d_{\bfX}}}  h\left(\wt \bx \right) c_{\bfY} \left(\frac 1\eta \wedge \|\bx \|^{a_\bfX} \right)^{\alpha_\bfY/a_\bfY} \: \P (\bX \in d \bx).
\end{equation*}
Now monotone convergence yields
\begin{equation*}
	\lim_{\eta \to 0}\lim_{t \to \infty} \frac{J_2(t)}{\P (\|\bfY\|^{a_\bfY} >t   )} = c_{\bfY}  \int_{\R^{d_{\bfX}}}  h\left(\wt \bx \right)  \|\bx \|^{a_\bfX\alpha_\bfY/a_\bfY} \: \P (\bX \in d \bx) =   c_{\bfY}\E h(\bfX).
\end{equation*}

\end{proof}

\section{Applications}\setcounter{equation}{0}\label{sec:RDE}
\subsection{Products of \regvary\ random matrices}\label{subsec:prod}
In what follows, we consider an iid \seq\ of 
$d\times d$ random matrices $(\bfA_i)$ and we assume that a generic element 
$\bfA\in\RV(\alpha,\mu^\bfA)$. We apply Theorem~\ref{thm:main1} to the 
\fct\ $\psi(\bfx,\bfy)= \bfx\cdot\bfy$.  

Next we formulate our findings for a general product 
$ {\mathbf \Pi}_n=\bfA_1\cdots\bfA_n$, $n\ge 1$. Here and in what follows,
we also use the notation
\beao
\mathbf \Pi_{i,j}=\left\{\barr{ll} \dprod_{s=i}^j \bfA_s\,,& i\le j\,,\\
{\rm Id}_d\,,&i>j\,, \earr\right.
\eeao
where ${\rm Id}_d$ is the $d\times d$ identity matrix.

\subsubsection{The case of non-equivalent tails} 
We first state the results in the case $\P(\|\mathbf\Pi_n\| >t) = o(\P(\|\mathbf\Pi_{n+1}\| >t))$ for all $n$. The complementary case is treated in Section~\ref{subsec:equiv}.

\bco\label{cor:product} Consider an iid \seq\ $(\bfA_i)$ of $d\times d$
matrices with $\bfA\in \RV(\alpha,\mu^\bfA)$. Assume that
\begin{equation}\label{eq:4:heavy}
	\dfrac{\P(\| \bfA \| >t)}{\P(\| \bfA_1 \| \cdot \| \bfA_2\| >t)} \to 0\,,\qquad t\to\infty\,.
\end{equation} 
Then for $n \geq 1$
\beam
\dfrac{\P\big(\|\mathbf\Pi_n\|>t\big)}{\P\big(\|\bfA_1\|\cdots \|\bfA_n\|>t\big)}\to \E\big[\|\bfTh_{\bfA_1}\cdots \bfTh_{\bfA_n}\|^\alpha\big]\,,\qquad t\to\infty\,.\label{eq:ko1}
\eeam 
If $\P(\|\bfTh_{\bfA_1}\cdots \bfTh_{\bfA_n}\| >0)>0$ then ${\mathbf \Pi}_n$ is \regvary\ and, as $t\to\infty$,
\beam
\lefteqn{\P\Big(\dfrac{\mathbf \Pi_n}{\|\mathbf\Pi_n\|}\in\cdot\,\Big|\,\|\mathbf \Pi_n\|>t\Big)}\nonumber\\ &\stw& \P\big(\bfTh_{\mathbf \Pi_n}\in\cdot\big) \label{eq:4:ind1}
=
\E\Big[  \dfrac{\|\bfTh_{\bfA_1}\cdots \bfTh_{\bfA_{n}}\|^\alpha}
{\E\big[ \|\bfTh_{\bfA_1}\cdots\bfTh_{\bfA_{n}}\|^\alpha\big]}
\1\Big( \dfrac{\bfTh_{\bfA_1} \cdots \bfTh_{\bfA_{n}}}{\|\bfTh_{\bfA_1}\cdots \bfTh_{\bfA_{n}}\|}\in \cdot\Big)\Big]\,.
\eeam 
In particular, if $\bfA$ is orthogonal,
\beao
\bfTh_{\mathbf \Pi_n}\stackrel{d}{=} \bfTh_{\bfA_1}\cdots \bfTh_{\bfA_n}\,.
\eeao
\eco
\bre In view of Lemma~\ref{lem:productregvaruniv}(2),
\eqref{eq:4:heavy} is satisfied if $\E[\|\bfA\|^\alpha] =\infty$.
\ere
\begin{proof}
We proceed by induction.  We will prove that for each $n$, \eqref{eq:4:ind1}, \eqref{eq:ko1} and
\begin{align}
& \P(\| \bfA_1 \| >t) + \P( \| \Pi_{2,n+1} \| >t) = o (\P (\|\bfA_1 \|\cdot \| \Pi_{2,n+1} \|>t))\label{eq:4:ind2}.
\end{align}
hold.
\par
We start with $n=2$. In view of~\eqref{eq:4:heavy} by Theorem~\ref{thm:main1}, 
\beao
\dfrac{\P\big(t^{-1} \bfA_1 \bfA_2\in \cdot\big)}{\P(\|\bfA_1\|\,\|\bfA_2\|>t)}
\stv   \E\big[\mu^\bfA(\{\bfx: \bfx\bfTh_\bfA\in\cdot\})\big]
\,.
\eeao
In particular, 
\beao 
\dfrac{\P\big(\|\bfA_1 \bfA_2\|>t\big)}{\P(\|\bfA_1\|\,\|\bfA_2\|>t)}
&\to& \E\big[\mu^\bfA(\{\bfx: \|\bfx \bfTh_\bfA\|>1\})\\
&=&\E\big[\mu^\bfA(\{\bfx: \|\wt \bfx \bfTh_\bfA\|>1/\|\bfx\|\})\\
&=&\int_{1}^\infty \P\big( r\,\|\bfTh_{\bfA_1}\bfTh_{\bfA_2}\|>1\big) d(-r^{-\alpha})\\
&=& \P\big( Y\,\|\bfTh_{\bfA_1}\bfTh_{\bfA_2}\|>1\big) = \E \big[\|\bfTh_{\bfA_1}\bfTh_{\bfA_2}\|^\alpha\big]\,,
\eeao
where $Y$ has a Pareto \ds , $\P(Y>r)=r^{-\alpha}$, $r>1$, independent
of the iid \rv s $\bfTh_{\bfA_1}$, $\bfTh_{\bfA_2}$. This proves~\eqref{eq:ko1} for $n=2$. 
Hence 
\beam\label{eq:feb10}
\dfrac{\P\big(t^{-1} \bfA_1 \bfA_2\in \cdot\big)}{\P(\|\bfA_1 \bfA_2\|>t)}\stv\dfrac{\E\big[\mu^\bfA(\{\bfx: \bfx\bfTh_\bfA\in\cdot\})\big]
 }{ \E \big[\|\bfTh_{\bfA_1}\bfTh_{\bfA_2}\|^\alpha\big]}\,.
\eeam
We conclude from \eqref{eq:feb10} that \eqref{eq:4:ind1} indeed holds for $n=2$ since
\beao
\P\Big(\dfrac{\bfA_1 \bfA_2}{\|\bfA_1\bfA_2\|}\in \cdot \,\Big|\,\|\bfA_1 \bfA_2\|>t\Big)&\stw&\dfrac{
\E\big[\mu^\bfA(\{\bfx: \dfrac{ \wt \bfx\bfTh_\bfA}{\|\wt \bfx\bfTh_\bfA\|}\in\cdot\,,\|\wt \bfx\bfTh_\bfA\|>1/\|\bfx\|\})\big]
}{  \E \big[\|\bfTh_{\bfA_1}\bfTh_{\bfA_2}\|^\alpha\big]     }\\
&=&\dfrac{\P\Big( Y \|\bfTh_{\bfA_1}\bfTh_{\bfA_2}\|\,\1\Big(
\dfrac{\bfTh_{\bfA_1}\bfTh_{\bfA_2}}{\|\bfTh_{\bfA_1}\bfTh_{\bfA_2}\|}\in\cdot\Big)
\Big)}{\E \big[\|\bfTh_{\bfA_1}\bfTh_{\bfA_2}\|^\alpha\big]}\\
&=& \E\Big[ \dfrac{\|\bfTh_{\bfA_1}\bfTh_{\bfA_2}\|^\alpha}{\E \big[\|\bfTh_{\bfA_1}\bfTh_{\bfA_2}\|^\alpha\big]} \,\1\Big(
\dfrac{\bfTh_{\bfA_1}\bfTh_{\bfA_2}}{\|\bfTh_{\bfA_1}\bfTh_{\bfA_2}\|}\in\cdot\Big)\Big]\\
&=&\P\big(\bfTh_{\bfA_1\bfA_2}\in\cdot\big)\,.
\eeao 
\par
To prove~\eqref{eq:4:ind2} for $n=2$  we note that we have already established
$$
	\P(  \| \mathbf\Pi_{2,3} \| >t) \sim  \E \big[\|\bfTh_{\bfA_1}\bfTh_{\bfA_2}\|^\alpha\big] \P( \| \bfA_1 \| \cdot \| \bfA_2 \| >t)  
$$ 
which, in combination with~\eqref{eq:4:heavy}, constitutes that for any $M>0$ there exists $t_0$ sufficiently large such that
\begin{equation*}
	\P( \| \mathbf\Pi_{2,3} \| >t) \geq M \,\P(\| \bfA_1 \|>t), \qquad t >t_0.
\end{equation*}
Take $\eta =  t_0^{-1}$.  We observe as $t\to\infty$ that
\begin{align*}
	\P (\|\bfA_1 \|\cdot \| \mathbf \Pi_{2,3} \|>t) & \geq \P (\|\bfA_1 \|\cdot \| \mathbf\Pi_{2,3} \|>t, \: \|\bfA_1\| \leq \eta t) \\
			& \geq M\, \P (\|\bfA_1 \|\cdot \| \bfA_2 \|>t, \: \|\bfA_1\| \leq \eta t) \\
			& \geq M\, \big( \P (\|\bfA_1 \|\cdot \| \bfA_2 \|>t) - \P (\|\bfA_1\| > \eta t)\big)\\
			& = M \,\P (\|\bfA_1 \|\cdot \| \bfA_2 \|>t)( 1+o(1)) \\
			& \geq M\, \P(\|\bfA_2\|>1)\P (\| \bfA_1 \|>t)( 1+o(1))\,. 
\end{align*}
The last two lines yield
\begin{equation*}
	\liminf_{t \to \infty} \frac{\P (\|\bfA_1 \|\cdot \| \mathbf \Pi_{2,3} \|>t)} {\P ( \| \mathbf\Pi_{2,3} \|>t)} \geq  \frac{M}{  \E \big[\|\bfTh_{\bfA_1}\bfTh_{\bfA_2}\|^\alpha\big]}
\end{equation*} 
and
\begin{equation*}
	\liminf_{t \to \infty} \frac{\P (\|\bfA_1 \|\cdot \| \mathbf \Pi_{2,3} \|>t)} {\P ( \| \bfA \|>t)} \geq  M\P(\|\bfA\|>1)\,,
	\end{equation*} 
	respectively. This proves \eqref{eq:4:ind2} for $n=2$ and finishes the proof 
of the corollary for $n=2$.
\par
Now suppose that it holds $n=k$ for some $k\ge 2$.
Since~\eqref{eq:4:ind2} holds for $n=k$ the balance conditions 
\beao
c_{\mathbf\Pi_{2,k+1}}&=& \lim_{t\to\infty} \dfrac{\P\big(\|\mathbf \Pi_{k}\|>t\big) }
{\P\big(\|\bfA_1\|\,\|\mathbf\Pi_{2,k+1}\|>t\big)}=0\,,\\
c_{\bfA_1}&=& \lim_{t\to\infty} \dfrac{\P\big(\|\bfA_1\|>t\big) }
{\P\big(\|\bfA_1\|\,\|\mathbf\Pi_{2,k+1}\|>t\big)}=0\,
\eeao
are satisfied. 
An application of Theorem~\ref{thm:main1} yields 
\beao
\dfrac{\P\big(t^{-1} \bfA_1 \mathbf\Pi_{2,k+1}\in \cdot\big)}{\P(\|\bfA_1\|\,\|\mathbf\Pi_{2,k+1}\|>t)}
\stv  \E\big[\mu^\bfA(\{\bfx: \bfx\bfTh_{\mathbf\Pi_k}\in\cdot\})\big]\,.
\eeao
An immediate con\seq \ is  
\beao
\dfrac{\P\big(\|\mathbf \Pi_{k+1}\|>t\big)}{\P(\|\bfA_1\|\,\|\bfA_2\cdots\bfA_{k+1}\|>t)}&\to&  
\E\big[\mu^\bfA(\{\bfx: \|\bfx\bfTh_{\mathbf\Pi_k}\|>1\})\big] 
=\P\big(Y \|\bfTh_{\bfA_1}\bfTh_{\mathbf\Pi_{2,k+1}}\|>1\big)\\
&=& \E\Big[ \|\bfTh_{\bfA_1}\bfTh_{\mathbf\Pi_{2,k+1}}\|^\alpha\Big]
=\dfrac{\E\big[ 
\|
\bfTh_{\bfA_1}\cdots\bfTh_{\bfA_{k+1}}
\|^\alpha\big]}{\E\big[ \|\bfTh_{\bfA_1}\cdots\bfTh_{\bfA_k}\|^\alpha\big]}\,,
\eeao
where the Pareto \rv\ $Y$, $\bfTh_{\bfA_1}$ and $\bfTh_{\mathbf\Pi_{2,k+1}}$ are 
independent. Here we also used the induction assumption on the \ds\
of $\mathbf \Pi_{k}$. Therefore 
\beao
\P\Big(\dfrac{\mathbf \Pi_{k+1}}{\|\mathbf\Pi_{k+1}\|}\in \cdot \,\Big|\,\|\mathbf\Pi_{k+1}\|>t\Big)
&\stw& \dfrac{\E\Big[\mu^\bfA\Big(\Big\{\bfx: 
\dfrac{\bfx\bfTh_{\mathbf \Pi_{k}}}{\|\bfx\bfTh_{\mathbf \Pi_k}\|}\in\cdot\,,\|\bfx\bfTh_{\mathbf\Pi_k}\|>1\Big\}\Big)\Big]}{ \E\big[ \|
\bfTh_{\bfA_1}\cdots\bfTh_{\bfA_{k+1}}
\|^\alpha\big]/\E\big[ \|\bfTh_{\bfA_1}\cdots\bfTh_{\bfA_k}\|^\alpha\big]}\\
&=&\dfrac{\P\Big( 
\dfrac{\bfTh_{\bfA_1} \bfTh_{\mathbf \Pi_{2,k+1}}}{
\|\bfTh_{\bfA_1}\bfTh_{\mathbf \Pi_{2,k+1}}\|}\in \cdot\,,Y\,\|\bfTh_{\bfA_1}\bfTh_{\mathbf\Pi_{2,k+1}}\|>1\Big)}{ \E\big[ 
\|
\bfTh_{\bfA_1}\cdots\bfTh_{\bfA_{k+1}}
\|^\alpha\big]/\E\big[ \|\bfTh_{\bfA_1}\cdots\bfTh_{\bfA_k}\|^\alpha\big]}\\
&=&
\dfrac{
\E\Big[\|\bfTh_{\bfA_1}\bfTh_{\mathbf\Pi_{2,k+1}}\|^\alpha\1\Big( \dfrac{\bfTh_{\bfA_1} \bfTh_{\mathbf \Pi_{2,k+1}}}
{\|\bfTh_{\bfA_1}\bfTh_{\mathbf \Pi_{2,k+1}}\|}\in \cdot\Big)\Big]}
{\E\big[ \|\bfTh_{\bfA_1}\cdots\bfTh_{\bfA_{k+1}}\|^\alpha\big]/
\E\big[ \|\bfTh_{\bfA_1}\cdots\bfTh_{\bfA_k}\|^\alpha\big]}
\\
&=&
\E\Big[{ \dfrac{ \|\bfTh_{\bfA_1}\cdots \bfTh_{\bfA_{k+1}}\|^\alpha}
{\E\big[ \|\bfTh_{\bfA_1}\cdots\bfTh_{\bfA_{k+1}}\|^\alpha\big]}}
\1\Big( \dfrac{\bfTh_{\bfA_1} \cdots \bfTh_{\bfA_{k+1}}}
{\|\bfTh_{\bfA_1}\cdots \bfTh_{\bfA_{k+1}}\|}\in \cdot\Big)\Big]
\,.
\eeao
This proves \eqref{eq:4:ind1} for $n=k+1$.  Finally, 
 we turn to \eqref{eq:ko1} for $n=k+1$:
\beao
\dfrac{\P\big(\|\mathbf\Pi_{ k+1}\|>t \big)}{\P\big(\|\bfA_1\|\cdots \|\bfA_{k+1}\|>t \big)}
&=&\dfrac{\P\big(\|\mathbf\Pi_{ k+1}\|>t \big)}{\P\big(\|\bfA_1\|\,\|\mathbf\Pi_{2,k+1}\|>t\big)}\dfrac{\P\big(\|\bfA_1\|\,\|\mathbf\Pi_{2,k+1}\|>t\big)}
{\P\big(\|\bfA_1\|\,\big(\|\bfA_2\|\cdots\|\bfA_{k+1}\|\big)>t\big)}\\
&\sim &\dfrac{\E\big[ 
\|
\bfTh_{\bfA_1}\cdots\bfTh_{\bfA_{k+1}}
\|^\alpha\big]}{\E\big[ \|\bfTh_{\bfA_1}\cdots\bfTh_{\bfA_k}\|^\alpha\big]}
\;\E\big[ \|\bfTh_{\bfA_1}\cdots\bfTh_{ \bfA_k}\|^\alpha\big]\\
&=&\E\big[ 
\|
\bfTh_{\bfA_1}\cdots\bfTh_{\bfA_{k+1}}
\|^\alpha\big]\,.
\eeao
In the last step we used the induction assumption leading to 
tail equivalence of $\|\bfA_2\|\,\|\bfA_3\cdots \bfA_{k+1}\|$
and $\|\bfA_2\,\bfA_3\cdots\bfA_{k+1}\|$
with factor $\E[ \|\bfTh_{\bfA_1}\cdots \bfTh_{\bfA_k}\|^\alpha]$.
To finish the proof we argue in favor of~\eqref{eq:4:ind2} for $n=k+1$ in the same fashion as we did that for $n=2$. More precisely, we have shown that
$$
	\P(  \| \mathbf\Pi_{k} \| >t) \sim  \E \big[\|\bfTh_{\bfA_1}\ldots \bfTh_{\bfA_k}\|^\alpha\big] \P( \| \bfA_1 \| \cdots \| \bfA_k \| >t)  
$$ 
which, in combination with~\eqref{eq:4:ind2} for $n=k$, gives 
$\P(\|\mathbf\Pi_k\| >t) = o (\| \mathbf \Pi_{k+1}\|>t))$. 
Consequently for any $M>0$  there exists $t_0$ sufficiently large such  that
\begin{equation*}
	\P( \| \mathbf\Pi_{k+1} \| >t) \geq M\, \P(\| \mathbf\Pi_k \|>t), \qquad t >t_0.
\end{equation*}
On the other hand, $  \P (\|\bfA_1\| > t) = o(\P (\|\bfA_1 \|\cdot \| \mathbf\Pi_{2, k+1} \|>t))$ and
\begin{equation*}
	\P(\| \mathbf\Pi_{k+1}\| > t) \sim  c_0\P(\|\bfA_1\| \cdot \| \mathbf\Pi_{2,k+1} \| >t), \qquad c_0 =\dfrac{\E\big[ 
\|
\bfTh_{\bfA_1}\cdots\bfTh_{\bfA_{k+1}}
\|^\alpha\big]}{\E\big[ \|\bfTh_{\bfA_1}\cdots\bfTh_{\bfA_k}\|^\alpha\big]}.
\end{equation*}
Take $\eta =  t_0^{-1}$. We observe as $t\to\infty$ that
\begin{align*}
	\P (\|\bfA_1 \|\cdot \| \mathbf \Pi_{2,k+2} \|>t) & \geq \P (\|\bfA_1 \|\cdot \| \mathbf\Pi_{2,k+2} \|>t, \: \|\bfA_1\| \leq \eta t) \\
			& \geq M \,\P (\|\bfA_1 \|\cdot \| \mathbf \Pi_{2,k+1} \|>t, \: \|\bfA_1\| \leq \eta t) \\
			& \geq M\,\big( \P (\|\bfA_1 \|\cdot \| \mathbf\Pi_{2, k+1} \|>t) - \P (\|\bfA_1\| > \eta t)\big)\\
			& = M \,\P (\|\bfA_1 \|\cdot \| \mathbf\Pi_{2,k+1} \|>t)( 1+o(1)) \\
			& = c_0\,M \,\P (\| \mathbf\Pi_{k+1} \|>t)( 1+o(1))\,. 
\end{align*}
This proves $ \P ( \| \mathbf \Pi_{k+1} \|>t) = o(\P (\|\bfA_1 \|\cdot \| \mathbf \Pi_{2,k+2} \|>t) )$ and finishes the proof of the corollary.
\end{proof}
\subsubsection{The case of tail-equivalent tails}\label{subsec:equiv}
We also assume condition \eqref{eq:mm}
which turns into
\beam\label{eq:mm1}
\lim_{M\to\infty}\limsup_{t\to\infty}\dfrac{\P(\|\bfA_1\|\,\|\bfA_2\|>t\,,M<\|\bfA_1\|\le t/M)}{\P(\|\bfA\|>x)}=0
\eeam
which is equivalent to
\begin{equation*}
	\frac{\P(\|\bfA_1\|>t)}{ \P(\|\bfA_1\|\cdot \|\bfA_2\|>t)}\to c_{\bfA}= \frac{1}{2 \,\E [\|\bfA\|^\alpha]}.
\end{equation*}
 An appeal to the following corollary shows that this condition causes
 tail equivalence of all $\mathbf\Pi_n$.
\bco
Consider an iid \seq\ $(\bfA_i)$ of 
$d\times d$ matrices such that $\bfA \in \RV(\alpha, \mu^\bfA)$ and 
\eqref{eq:mm1} holds.
Then for any $n\geq 2$,
\begin{equation}\label{eq:4:cos}
\frac{\P \left(  \|\mathbf \Pi_n\| >t \right)}{ \P (\| \bfA \| >t)}\to  \sum_{k=1}^{n} \E\big[ \|\mathbf\Pi_{k-1}\bfTh_{\bfA_{k}} \mathbf\Pi_{k+1,n}\|^\alpha\big]\,,\qquad t\to\infty\,.
\end{equation}
Additionally, if  $\P(\|\mathbf\Pi_{k-1}\bfTh_{\bfA_{k}} \mathbf\Pi_{k+1,n}\| >0)>0$ for some $k\leq n$ then ${\mathbf \Pi}_n$ is \regvary\ and as $t\to\infty$,
\beao
\lefteqn{\P\Big(\dfrac{\mathbf \Pi_n}{\|\mathbf\Pi_n\|}\in\cdot\,\Big|\,\|\mathbf \Pi_n\|>t\Big) \stw \P\big(\bfTh_{\mathbf \Pi_n}\in\cdot\big)}\\
&=& 
\sum_{k=1}^n p_k\,\E\Big[ \dfrac{\|\mathbf\Pi_{k-1}\bfTh_{\bfA_{k}} \mathbf\Pi_{k+1,n}\|^\alpha}{\E\big[\|\mathbf\Pi_{k-1}\bfTh_{\bfA_{k}} \mathbf\Pi_{k+1,n}\|^\alpha\big]}\1\Big( \dfrac{\mathbf\Pi_{k-1}\bfTh_{\bfA_{k}} \mathbf\Pi_{k+1,n}}{\|\mathbf\Pi_{k-1}\bfTh_{\bfA_{k}} \mathbf\Pi_{k+1,n}\|}\in \cdot\Big)\Big]
\eeao
where 
\beao
p_k= \dfrac{\E\big[\|\mathbf\Pi_{k-1}\bfTh_{\bfA_{k}} \mathbf\Pi_{k+1,n}\|^\alpha\big]}{\sum_{k=1}^n\E\big[ \|\mathbf\Pi_{k-1}\bfTh_{\bfA_{k}} \mathbf\Pi_{k+1,n}\|^\alpha\big]}\,,\qquad k=1,\ldots,n\,. 
\eeao 
\eco
\begin{proof}
	We proceed by induction. We will prove~\eqref{eq:4:cos} and 
	\beao\mu^{\mathbf \Pi_n}(\cdot ) &=&\dfrac{
	 \sum_{k=1}^n \E\Big[ \mu^\bfA \left( \mathbf a  \: : \: \mathbf\Pi_{k-1}\mathbf a\mathbf\Pi_{k+1,n} \in \cdot\right)\Big]}
{\sum_{k=1}^n\E\big[ \|\mathbf\Pi_{k-1}\bfTh_{\bfA_{k}} \mathbf\Pi_{k+1,n}\|^\alpha\big]}\,.
	\eeao
\par
For $n=2$,  Theorem~\ref{thm:main1} yields
\beao
\dfrac{\P\big(t^{-1}\bfA_1\bfA_2\in\cdot\big)}{\P(\|\bfA\|>t)}
&=&\dfrac{\P\big(t^{-1}\bfA_1\bfA_2\in\cdot\big)}
{\P\big(\|\bfA_1\|\,\|\bfA_2\|>t\big)}\dfrac{\P\big(\|\bfA_1\|\,\|\bfA_2\|>t\big)}{\P(\|\bfA\|>t)}\\
&\stv& \E\big[\mu^\bfA(\{\bfx: \bfx\,\bfA\in\cdot\})\big]+
\E\big[\mu^\bfA(\{\bfx: \bfA\,\bfx\in\cdot\})\big]\\
&=& \int_0^{\infty}\alpha r^{-\alpha-1}\big(\P\big(r\bfTh_{\bfA_1}\bfA_2\in\cdot\big)+\P\big(r\bfA_1\bfTh_{\bfA_2}
\in \cdot\big)\big) \, dr\,.
\eeao
In particular, for a Pareto \rv\ $Y$ independent of $\bfA_1,\bfA_2$
and $\bfTh_{\bfA_1},\bfTh_{\bfA_2}$,
\beao
\dfrac{\P\big(\|\bfA_1\bfA_2\|>t\big)}{\P\big(\|\bfA\|>t\big)}
&\to&  \P\big(Y\|\bfTh_{\bfA_1}\bfA_2\|>1\big)+\P\big(Y\|\bfA_1\bfTh_{\bfA_2}\|
>1\big)\\
&=& \E\big[\|\bfTh_{\bfA_1}\bfA_2\|^\alpha+\|\bfA_1\bfTh_{\bfA_2}\|^\alpha \big]\,.
\eeao
We also have 
\beao
\lefteqn{\P\Big(\dfrac{\bfA_1\bfA_2}{\|\bfA_1\bfA_2\|}\in\cdot \,\Big|\,\|\bfA_1\bfA_2\|>t\Big)}\\
&\stw& \dfrac{ \P\Big(\dfrac{\bfTh_{\bfA_1}\bfA_2}{\|\bfTh_{\bfA_1}\bfA_2\|}\in\cdot\,,Y\,\|\bfTh_{\bfA_1}\bfA_2\|>1\Big)+\P\Big(\dfrac{\bfA_1\bfTh_{\bfA_2}}{\|\bfA_1\bfTh_{\bfA_2}\|}
\in \cdot\,,  Y\,\|\bfA_1\bfTh_{\bfA_2}\|>1\Big)}{\E\big[\|\bfTh_{\bfA_1}\bfA_2\|^\alpha+\|\bfA_1\bfTh_{\bfA_2}\|^\alpha\big]}\\
&=& \dfrac{ \E \Big[\|\bfTh_{\bfA_1}\bfA_2\|^\alpha\1\Big(\dfrac{\bfTh_{\bfA_1}\bfA_2}{\|\bfTh_{\bfA_1}\bfA_2\|}
\in \cdot \Big) + \E \Big[\|\bfA_1\bfTh_{\bfA_2}\|^\alpha\1\Big(\dfrac{\bfA_1\bfTh_{\bfA_2}}{\|\bfA_1\bfTh_{\bfA_2}\|}
\in \cdot \Big)
\Big]}{\E\big[\|\bfTh_{\bfA_1}\bfA_2\|^\alpha+\|\bfA_1\bfTh_{\bfA_2}\|^\alpha\big]}\,.
\eeao
\par
Now suppose that our claim  holds for some $n\geq 2$. Put $\tilde c_n =  \sum_{k=1}^{n} \E\big[ \|\mathbf\Pi_{k-1}\bfTh_{\bfA_{k}} \mathbf\Pi_{k+1,n}\|^\alpha\big]$.
Since $\|\bfA\|$ satisfies~\eqref{eq:mm1}  and $\P(\| \mathbf \Pi_n \| > t) \sim \tilde c_n\, \P(\|\bfA\|>t)$ we infer that
\beao
	\dfrac{\P(\|\bfA \| >t)}{\P(\|\bfA\|\cdot \|\mathbf\Pi_{2,n+1} \| >t)  } &\to & c_{n, \bfA} = \frac{1}{\E [\|\mathbf \Pi_{2, n+1} \|^\alpha]  + \tilde c_n\, \E [ |\bfA\|^\alpha]}\,,\\
	\dfrac{\P(\|\mathbf\Pi_{2,n+1}\| >t)}{\P(\|\bfA\|\cdot \|\,\mathbf\Pi_{2,n+1} \| >t)  } &\to&  c_{n, \mathbf \Pi} = \dfrac{\tilde c_n}{\E[ \|\mathbf \Pi_{2, n+1} \|^\alpha]  + \tilde c_n\, \E[  \|\bfA\|^\alpha]}.
\eeao	
Theorem~\ref{thm:main1} yields
\beao
\dfrac{\P(t^{-1}\mathbf \Pi_{n+1} \in \cdot )}
{\P(\|\bfA\|\cdot \|\mathbf\Pi_{2,n+1} \| >t)} 
\to 
\dfrac{\E\Big[ \mu^\bfA\big(\big\{ \mathbf a \: : \: 
\mathbf a \,\mathbf\Pi_{2,n+1} \in \cdot \big\}\big) + \tilde c_n\,
\mu^{\mathbf\Pi_n} \big(\big\{{\mathbf \pi} \:: \: \bfA \,{\mathbf\pi} \in \cdot\big\} \big)\Big]  }
{\E \big[\|\mathbf \Pi_{2, n+1} \|^\alpha\big]  + \tilde c_n\, \E[ \|\bfA\|^\alpha]}\,.
\eeao
Consequently, by the induction hypothesis,
\beao
\dfrac{\P(t^{-1} \mathbf\Pi_{n+1} \in \cdot )}{\P( \|\mathbf\Pi_{n+1} \| >t) } 
& \stv& 
\dfrac{\E \Big[\mu^\bfA\big(\big\{ \mathbf a \: : \: \mathbf a\, 
\mathbf\Pi_{2,n+1} \in \cdot \big\}\big) + \tilde c_n\,
\mu^{\mathbf\Pi_n} \big(\big\{{\mathbf \pi} \:: \: \bfA\, {\mathbf\pi} \in \cdot 
\big\}\big)\Big]}{ \E\Big[ \mu^\bfA\big(\big\{ \mathbf a \: : \: 
\| \mathbf a \,\mathbf\Pi_{2,n+1}\|>1 \big\}\big) + \tilde c_n\,
\mu^{\mathbf\Pi_n} \big(\big\{{\mathbf \pi} \: \: \|\bfA \,{\mathbf\pi}\|>1\big\}\big)\Big]} \\
& =&\dfrac{\E \Big[\mu^\bfA\big(\big\{ \mathbf a \: : \: 
\mathbf a \,\mathbf\Pi_{2,n+1} \in \cdot\big\}\big) + 
 \sum_{k=1}^n\mu^{\bfA} \big(\big\{\mathbf a \: : \:  
\mathbf\Pi_{k}\,\mathbf a\,\mathbf\Pi_{k+2,n+1} \in \cdot\big\}\big )\Big]}
{ \E\Big[ \mu^\bfA\big(\big\{ \mathbf a \: : \: \| \mathbf a\, \mathbf\Pi_{2,n+1}\|>1\big\}\big ) + 
\sum_{k=1}^n\mu^{\bfA} \big(\big\{\mathbf a \:: \:  
 \|\mathbf\Pi_{k}\,\mathbf a\,\mathbf\Pi_{k+2,n+1}\| >1\big\}\big)\Big]} \\
	& =&\dfrac{ \sum_{k=1}^{n+1} \E\big[\mu^{\bfA} \big(\big\{\mathbf a \:: \: 
\mathbf\Pi_{k-1}\,\mathbf a\, \mathbf\Pi_{k+1,n+1} \in \cdot\big\}\big)\big]}{\sum_{k=1}^{n+1}\E\big[ \|\mathbf\Pi_{k-1}\,\bfTh_{\bfA_{k}}\,\mathbf\Pi_{k+1,n+1}\|^\alpha\big]}\,.
\eeao
With this at hand, the convergence 
\beao
\P\Big(\dfrac{\mathbf \Pi_n}{\|\mathbf\Pi_n\|}\in\cdot\,\Big|\,\|\mathbf \Pi_n\|>t\Big) \stw\P\big(\bfTh_{\mathbf \Pi_n}\in\cdot\big)
\eeao
follows.
\end{proof}

\subsection{Stochastic recurrence equations}\label{sec:RDE1}
We  turn to the \sre\ 
\beam\label{eq:22a}
\bfR_t=\bfA_t\,\bfR_{t-1}+\bfB_t\,,\qquad t\in\bbz\,,
\eeam
where $\big((\bfA_t,\bfB_t)\big)_{t\in\bbz}$ is an iid \seq\ with generic element
$(\bfA,\bfB)$, $\bfA$ is a $d\times d$ random matrix and $\bfB$  
an $\bbr^d$-valued random vector, possibly dependent  on each other. A solution $(\bfR_t)$
is causal if for every $t$, $\bfR_t$ is a \fct\ only of 
values $\big((\bfA_s,\bfB_s)\big)_{s\le t}$, and then 
it constitutes a \MC . If a stationary causal 
solution $(\bfR_t)$ with generic element $\bfR$ exists its marginal \ds\ satisfies the fixed point equation in law
	\beam\label{eq:22}
		\bR\eqd \bfA\,\bR+\bfB\,, 
	\eeam 
and $\bfR$ has the \rep\ in law
	\begin{equation}\label{eq:defR}
		\bR\eqd \sum_{k=0}^\infty \mathbf \Pi_k \bfB_{k+1}\,, \qquad\mbox{ with $\mathbf \Pi_k = \prod_{j=1}^k \bfA_j$}\,.
	\end{equation}
The latter infinite series converges under  conditions on the \ds\ of 
$(\bfA,\bfB)$,  for example $\E [\log \| \bfA \|] <0$ and 
$\E [\log_+\|\bfB\|] <\infty$. 
	Under some mild integrability and non-degeneracy assumptions~\eqref{eq:defR} is the unique solution to~\eqref{eq:22}. Here and in what follows,
we refer to the monograph
Buraczewski~et~al.~\cite{buraczewski:damek:mikosch:2016} for 
	details concerning the existence, uniqueness and other properties
of the solutions to~\eqref{eq:22a} and \eqref{eq:22}.
\par	
The equations \eqref{eq:22a} and \eqref{eq:22} have 
attracted a lot of attention since the seminal paper by 
Kesten \cite{kesten:1973} who proved that  
	$\bR$ has some \regvar\ property with tail index $\alpha>0$ given by
	\begin{equation*}
		\lim_{n \to \infty} (\E \| \mathbf \Pi_n\|^\alpha)^{1/n}=1.
	\end{equation*}
If $d=1$, the latter equation reads as $\E[|A|^\alpha]=1$.
	 In the Kesten setting, it is typically assumed 
that $\E [\|\bfB\|^\alpha]<\infty$ and $\E [\|\bfA\|^\alpha\log_+ \|\bfA\|]<\infty$,
implying the existence and uniqueness of the solution $(\bfR_t)$.
Under these and further mild conditions on the \ds\ of 
$(\bfA,\bfB)$ one has $\bfR\in \RV(\alpha,\mu^\bfR)$ and the tail \asy s 
\beao
\P(\|\bfR\|>t)\sim c_0\,t^{-\alpha}\qquad\mbox{ for some $c_0>0$}\,.
\eeao
Since $\E[\|\bfR\|^\alpha]=\infty$ we have 
$\P(\|\bfB\|>t)=o(\P(\|\bfR\|>t))$, and  elementary calculations 
	(Lemma~C.3.1 in Buraczewski et al. \cite{buraczewski:damek:mikosch:2016}) show that for $\mu^\bfR$-continuity sets $C$,
	\beao
		t^\alpha \P(t^{-1}\bR\in C )\sim t^{\alpha}\,\P(t^{-1}\bfA\,\bR\in C)\,,
	\eeao
and the multivariate Breiman result Lemma~C.3.1 in  \cite{buraczewski:damek:mikosch:2016} yields
 	\beao
		\dfrac{\P(t^{-1}\bfA\,\bR\in \cdot)}{\P(\|\bR\|>t)}\stv  \E \big[\mu^\bR\big(\{\bfx: \bfA\bfx \in \cdot\big)\}\big)\big]\,.
	\eeao
	Hence we have the identity
	\beao
		\mu^\bR(\cdot) = \E \big[\mu^\bR\big(\{\bfx: \bfA\bfx \in \cdot \}\big)\big]\,.
	\eeao
Using induction on the recursion \eqref{eq:22a} and similar arguments, we find that
\beao
		\mu^\bR(\cdot) = \E \big[\mu^\bR\big(\{\bfx: {\mathbf \Pi}_k\bfx \in \cdot \}\big)\big]\,,\qquad k\ge 1\,.
\eeao
This relation holds, in particular, if $\bfA$ is \regvary\ with index $\alpha$ but the additional 
moment condition $\E [\|\bfA\|^\alpha\log_+ \|\bfA\|]<\infty$ must be satisfied.
\par
Regular variation of $(\bfR_t)$ may also arise from \regvar\ of $\bfB$
under the alternative conditions 
\beam\label{eq:dysz}
\bfB \in  \RV \big( \alpha,\mu^\bfB\big),  \quad
\E[\|\bfA\|^{\alpha}]<1\,\quad\mbox{\,and 
	$\E[\|\bfA\|^{\alpha+\delta}]<\infty$ for some $\delta>0$}\,.
\eeam
Then $\bR$ is regularly varying with index $\alpha$ and
	\begin{equation*}
		\dfrac{\P\left( t^{-1} \bR \in \cdot \right)}{ \P(\| \bfB \| > t)} \stv \int \mu_\bfB ( \{ \by\: : \: \bz\by \in \cdot  \}) \: \nu_{\mathbf \Pi}(d \bz),
	\end{equation*}
	where $\nu_{\mathbf \Pi}(\cdot) = \sum_{k=0}^\infty \P (\mathbf{\Pi}_k \in \cdot)$ ia s measure on $\M_{d \times d}$; see Theorem~4.4.24 
	in \cite{buraczewski:damek:mikosch:2016}.
\par

For our purposes we will treat $(\bfA,\bfB)$ as a random element of ${\mathbb M}_{d\times d}\times \bbr^d$ equipped with the norm $\| (\bf a, \bf b) \| = \|\bf a \| + \|\bf b\|$, where  
$\|\bf a \|$ stands for the operator norm of the matrix $\bf a$ (with respect to the Euclidean distance) and   $\|\bf b\|$ is the Euclidean norm of the vector $\bf b$. 
We assume that the following set of conditions $(\bfC)$ 
on $(\bfA,\bfB)$ holds:
\begin{enumerate}
\item[\rm (C1)] A \regvar\ condition holds for some non-null Radon 
\ms\ $ \mu^{(\bfA,\bfB)}$ on ${\mathbb M}_{d\times d}\times \bbr^d$:
\beam\label{eq:p}
\dfrac{\P\big(t^{-1} (\bfA,\bfB)\in\cdot\big)}{\P(\|(\bfA, \bfB)\|>t)}\stv 
\mu^{(\bfA,\bfB)}(\cdot)\,,\qquad t\to\infty\,.
\eeam
\item[\rm (C2)] $X=\|(\bfA_1, \bfB_1)\|$ and $Y=\|(\bfA_2, \bfB_2)\|$ satisfy \eqref{eq:strange}.
\item[\rm (C3)] $\E[\|\bfA\|^{\alpha}]<1$ and $\mu^{(\bfA, \bfB)}\big(\big\{\mathbf{a}, \mathbf{b})\: : \: \|\mathbf{a}\|>1\big\}\big) >0$.
\end{enumerate}
\subsubsection*{Some comments}
\begin{itemize}
\item
To the best of our knowledge, except for some univariate cases 
treated in Damek and Dyszewski~\cite{damek:dyszewski:2018} and Kevei~\cite{kevei:2016},  
not much is known about \regvar\ of  $\bfR$ under \regvar\ of $\bfA$ and 
(C3). Then \eqref{eq:dysz} is violated since $\E[\|\bfA\|^{\alpha+\delta}]=\infty$ for any 
	$\delta>0$. 
\item
In view of Lemma~\ref{lem:productregvaruniv} condition (C2) implies
\beao
		\dfrac{\P( \|(\bfA_1, \bfB_1)\| \cdot \| (\bfA_2, \bfB_2)\|>t)}{\P(\|(\bfA, \bfB) \| >t)} \to 2\,\E \|(\bfA, \bfB)\|^{\alpha}\,.
\eeao
\end{itemize}
\par
The following result is a  multivariate counterpart of the 
results obtained in~Damek and Dyszewski~\cite{damek:dyszewski:2018}.
\bth\label{thm:main2}
Assume {\rm $(\bfC)$}. Then $\bfR$ given in \eqref{eq:defR} satisfies  
\beao
\dfrac{\P(t^{-1}\bR\in \cdot)}{\P(\|(\bfA, \bfB\|>t)}\stv  
\nu(\cdot) = \sum_{n=0}^{\infty} \E \big[\mu^{(\bfA, \bfB)} ( \{ (\ba, \bb) \: : \: \mathbf \Pi_n (\ba \,\bR_0 + \bb) \in \cdot \})\big]\,.
\eeao
		In particular, if the measure $\nu$ on $\R^d_{\bf0}$ is non-null
 then $\bR \in  \RV \big( \alpha,\mu^\bR\big)$ with
		\begin{equation*}
			\mu^{\bR}(\cdot) = \dfrac{\nu(\cdot)} {\nu\big(\big\{ \br\: : \: \|\br\|>1 \big\}\big)}\, .
		\end{equation*}
\ethe
The remainder of this section is devoted to the proof of the theorem.
A main step in the proof is provided by the following lemma.
\ble\label{lem:iter}
Assume that the $\R^d$-valued  random vector  
$\bfX\in {\rm RV}(\alpha,\mu^\bfX)$ is   independent of $(\bfA, \bfB)$ which
satisfies~{\rm \bf (C)} and there is a positive constant $d_{\bfX}$ such that
\beam\label{eq:dx}
\dfrac{\P( \|\bfX\|>t  )}{\P(\|(\bfA,\bfB)\|>t)} \to d_\bfX\,,\qquad t\to\infty\,.
\eeam
Then as $t\to\infty$,
\beao
			\dfrac{\P( \| \bfA\bfX + \bfB \| >t)}{\P(\|(\bfA,\bfB)\|>t)} &\to& \E \big[\mu^{(\bfA, \bfB)} (\{ (\ba,\bb) \: : \: \|\ba\bfX+\bb\|>1\})\big] 
				+  d_\bfX\,\E \big[ \| \bfA \|^\alpha \big]=:C_0\,,\\
\dfrac{\P(t^{-1}(\bfA\bfX + \bfB) \in \cdot)}{\P(\|\bfA \bX+\bfB\|>t)}&\stv& 
C_0^{-1}\,\E \big[  \mu^{(\bfA, \bfB)} (\{ (\ba,\bb) \: : \: \ba\bfX+\bb \in \cdot\})\big] 
				+  d_\bfX\,\E \big[\mu^{\bfX} (\{ \bx \: : \: \bfA\bx\ \in \cdot \})\big]\,.
	\eeao
	\ele

	\begin{proof}[Proof of Lemma~\ref{lem:iter}]
Write $\1_d =(1,\ldots,1)^\top\in \R^d$, ${\rm Id}_d$ and  $\diag(\bb)$, $\bb \in \R^d$, in $\mathbb M_{d\times d}$ for the identity matrix and 
the diagonal matrix whose 
consecutive diagonal entries are the consecutive components of $\bb$, respectively.
We write
		\begin{equation*}
			\wh \bfX = \left( \begin{array}{c} \bfX \\ \1_d\end{array} \right) \in \R^{2d} \quad \mbox{and} \quad 
			\wh \bfA = \left( \begin{array}{cc} \bfA & \diag(\bfB) \\  \mathbf 0 & {\rm Id}_d\end{array} \right) \in \mathbb M_{2d\times 2d},
		\end{equation*}
Then $\wh \bfX$ and $\wh\bfA$ are both regularly varying. Indeed, for 
$\wh \bfX$  we have
		\beao
\dfrac{\P(t^{-1}\wh \bfX\in \cdot)}{\P(\|\wh \bfX\|>t)}
&\sim&\dfrac{\P(t^{-1}\wh \bfX\in \cdot)}{\P(\|\bfX\|>t)}\\
&\stv&\mu^{\wh \bfX}(\cdot)
 = \mu^{\bfX} \Big( \Big\{ \bx\in\bbr^d_{\bf0} \: : \: 
\Big( \barr{c} \bfx \\ {\bf0}\earr \Big) \in \cdot \Big\}\Big)\,.
		\eeao
For $\wh \bfA$, choosing the operator norm $\|\cdot\|$, we have
\beao
\dfrac{\P(t^{-1}\wh \bfA\in \cdot)}{\P(\|\wh \bfA\|>t)}\stv 
\mu^{\wh \bfA}(\cdot) = 
\dfrac{ \mu^{(\bfA, \bfB)} \Big( \Big\{ (\ba, \bb) \: : \:  \Big( \begin{array}{cc} \ba & \diag(\bb) \\  \mathbf 0 & \mathbf 0\end{array} \Big) \in \cdot\Big\}\Big) }
			{ \mu^{(\bfA, \bfB)} \Big(\{ (\ba, \bb) \: : \:  \|\ba\| \vee \| \diag(\bb)\| >1 \} \Big)}.
\eeao
We intend to use the fact that
\begin{equation*}
\wh \bfA\, \wh\bfX = \Big( \begin{array}{c} \bfA\bfX + \bfB \\ \1_d\end{array} \Big)
\end{equation*}
		in combination with Theorem~\ref{thm:main1} to prove 
the claim. In view of the tail equivalence condition \eqref{eq:dx} we have
\beao
c_{ \wh \bfA}& =& \lim_{t \to \infty } \dfrac{\P( \| \wh \bfA \| > t )}{\P( \| \wh \bfA \| \cdot \| \wh \bfX \| > t )} = \dfrac{1}{\E [\| \wh \bfX \|^\alpha] + d_\bfX \E [\| \wh \bfA\|^\alpha] }\\
c_{ \wh \bfX} &=& \lim_{t \to \infty } \dfrac{\P( \| \wh \bfX \| > t )}{\P( \| \wh \bfA \| \cdot \| \wh \bfX \| > t )} = \dfrac{d_\bfX}{\E[ \| \wh \bfX \|^\alpha] + d_\bfX \,\E [\| \wh \bfA\|^\alpha] }.
\eeao
Therefore Theorem~\ref{thm:main1} yields
\beao
\lefteqn{\dfrac{\P(t^{-1}(\bfA\bfX + \bfB) \in \cdot)}{\P(\|(\bfA,\bfB)\|>t)}}\\  &\stv  & c_{\wh\bfA}\, \E \big[\mu^{\wh\bfA}\big(\big\{\wh{\mathbf{a}} \: : \: \wh{\mathbf{a}} \wh\bfX \in \cdot\big\}\big )\big] 
				+  c_{\wh\bfX}\, \E\big[ \mu^{\wh\bfX}\big(\big\{\wh{\mathbf{x}} \: : \: \wh{\mathbf{A}}\, \wh{\mathbf{x}} \in \cdot\big\}\big )\big] \\
				 & =& c_{\wh\bfA}\, \E \big[\mu^{(\bfA, \bfB)}\big(
\big\{(\mathbf{a}, \mathbf{b}) \: : \: \mathbf{a} \,\wh\bfX + \mathbf{b} \in \cdot\big\}\big )\big] +  c_{\wh\bfX}\, \E\big[ \mu^{\bfX}\big(\big\{\mathbf{x} \: : \: \mathbf{A}\, \mathbf{x} \in \cdot\big\}\big )\big] \\
&=& \dfrac{\E\big[ \mu^{(\bfA, \bfB)}\big(\big\{(\mathbf{a}, \mathbf{b}) \: : \: \mathbf{a} \wh\bfX + \mathbf{b} \in \cdot\big\}\big)\big] 
				+  d_{\bfX}\, \E\big[ \mu^{\bfX}\big(\big\{\mathbf{x} \: : \: \mathbf{A}\, \mathbf{x} \in \cdot \big\}\big)\big]}{\E [\| \wh \bfX \|^\alpha] + d_\bfX\, \E [\| \wh \bfA\|^\alpha]}\,,
\eeao
		which implies both claims.
	
\end{proof}	
Consider the Markov chain $(\bR_n^0)_{n \geq 0}$ given by the
recursion \eqref{eq:22a}  with $\bR_0^0 =0$. Then
\beao
\bR_n^0 \eqd \sum_{k=0}^{n-1} \mathbf \Pi_k \bfB_{k+1}\std \bR\,.
\eeao
By Lemma~\ref{lem:iter},
	\begin{equation*}
			\frac{\P( t^{-1} \bR_n^0  \in \cdot )}{\P(\|( \bfA, \bfB) \| > t )} \stv \nu_n(\cdot),
	\end{equation*}
and the sequence $(\nu_n)_{n \geq 0}$ of measures on $\R^d_{\bf0}$ 
satisfies the recursive relation 
	\beam\label{eq:receta}
\nu_{n+1}(\cdot) &=& \E\big[ \mu^{(\bfA, \bfB)} ( \{ (\ba,\bb) \: : \: \ba\bR_n^0+\bb \in \cdot \})\big] +  \E\big[ \nu_n (\{ \bx \: : \: \bfA\bx\ \in \cdot \})\big]\,,\quad n\ge 1\,,\nonumber\\
\nu_0&=&o\,.\nonumber\\
	\eeam
We have
\begin{equation*}
\Big\|\sum_{k=0}^\infty \mathbf \Pi_{k}\bfB_{k+1}\Big\|\le 	R = \sum_{k=0}^\infty \| \bfB_{k+1}\| \prod_{j=1}^k \| \bfA_j\|
\end{equation*}
A copy $\wt R$ of $R$ which is also independent of $(\bfA,\bfB)$ 
solves the equation
	\begin{equation*}
		\wt R \stackrel{d}{=} \|  \bfA \| \,\wt R + \| \bfB \|\,,
	\end{equation*}
and $\| \bR_n^0 \| \stackrel{d}{\leq} R$, $n\ge 0$, 
where for any non-negative \rv s $X,Y$, $ X\stackrel{d}{\leq}Y$ 
stands for stochastic domination, i.e., 
$\P(Y > t) \ge \P(X > t)$ for any $t>0$.
From the main result in Damek and Dyszewski 
\cite{damek:dyszewski:2018} (see Lemma \ref{lem:dd2018}) we also have
under {\bf \rm (C)},
\beam
\limsup_{t \to \infty}\frac{\P(\| \bR\| >t)}{\P(\| \bfA\| >t)}& \leq&  \limsup_{t \to \infty}\frac{\P( R >t)}{\P(\| \bfA\| >t)}  < \infty,\nonumber\\
\sup_n \E \big[\|\bR_n^0 \|^\alpha\big]& \leq &\E \big[R^\alpha]<\infty .\label{eq:d}
		\eeam
\ble\label{lem:last}
Assume {\bf \rm (C)}. Then 
\begin{equation*}
			\nu_n(\cdot) \stv \nu(\cdot) = \sum_{k=0}^{\infty} \E \mu^{(\bfA, \bfB)} ( \{  (\ba, \bb) \: : \: \mathbf \Pi_k (\ba \bR_0 + \bb) \in \cdot \} ),
\end{equation*}
where $\nu$ is a Radon measure on $\R^d_{\bf0}$.
	\ele	
\begin{proof}[Proof of Lemma~\ref{lem:last}] For $k\le n$
write $\mathbf \Pi^\downarrow_{n,k}= \bfA_n \bfA_{n-1} \cdots \bfA_k$.
We have by \eqref{eq:receta}, 
		\begin{equation*}
			\nu_n(\cdot) = \sum_{k=1}^{n} \E 
\big[\mu^{(\bfA, \bfB)} ( \{ (\ba, \bb) \: : \: \mathbf \Pi^\downarrow_{n,k+1} (\ba \bR_{k-1}^0 + \bb) \in \cdot \} )\big] =  \sum_{k=1}^{n} \eta_{n,k}(\cdot).
		\end{equation*}
Write
		\begin{equation*}
\nu(\cdot) = \sum_{k=0}^{\infty} \E \big[\mu^{(\bfA, \bfB)} ( \{  (\ba, \bb) \: : \: \mathbf \Pi_k (\ba \bR_0 + \bb) \in \cdot \} )\big] =  \sum_{k=0}^{\infty}  \eta_k(\cdot),
		\end{equation*}
We intend to show $\nu_n\stv\nu$ or, equivalently,
 $\int f d\nu_n\to \int fd\nu$ for any $f \in C_c^+(\R^d_{\bf0})$. Then there 
are $c,M>0$ such that $f$ vanishes  on $\{\bfx:\|\bfx\|>c\}$
and $f(\bfx)\le M<\infty$.
Our strategy is to use the following approximations: 
\begin{equation*}
	\int f d\nu_n \stackrel{\bf (1)}{\approx} \int f d \Big( \sum_{ n/2< k \leq n} \eta_{n,k}\Big) \stackrel{\bf (2)}
{\approx}  \int f d \Big( \sum_{ 0< k \leq n/2} \eta_{k}\Big) \stackrel{\bf (3)}{\approx} \int f d\nu.
\end{equation*}
In what follows, we will make these approximations precise.\\[2mm]
{\bf Approximations (1) and (3)}. For {\bf (1)}, we will show that
\beam
\lim_{n\to \infty}  \int f d \Big(\sum_{ k \leq n/2} \eta_{n,k}\Big) 
&=& \lim_{\nto}\sum_{k=1}^{ [ n/2]} \E\Big[  \dint  
f\big( \mathbf \Pi^\downarrow_{n,k+1} (\ba \bR_{k-1}^0 + \bb)\big)\,   
\mu^{(\bfA, \bfB)}(d (\ba,\bb))\Big]=0\,,\nonumber\\\label{eq:a}
\eeam
For any $c>0$ and $k\le [n/2]$ we have
\beao
\lefteqn{\E 
\big[\mu^{(\bfA, \bfB)} \big( \big\{ (\ba, \bb) \: : \: \|\mathbf \Pi^\downarrow_{n,k+1} (\ba \bR_{k-1}^0 + \bb)\|> c \big\} \big)\big]}\\
&\le &\E 
\big[\|\mathbf \Pi^\downarrow_{n,k+1}\|^\alpha\big]\,\E\big[\mu^{(\bfA, \bfB)} \big( \big\{ (\ba, \bb) \: : \: \| \ba \bR_{k-1}^0 + \bb\|> c \big\} \big)\big]\\
		 		 &\leq &  
(\E\big[ \| \bfA \|^\alpha\big])^{n-k}\,\Big( \E\big[ \mu^{(\bfA, \bfB)} \big(\big \{ (\ba, \bb) \: : \: \|\ba \bR_{k-1}^0\| > c/2\big\}\big)\big]+ 
\mu^{(\bfA, \bfB)} \big( \big\{ (\ba, \bb) \: : \: \|\bb\| > c/2\big\}\big)
\Big)\\
&\le &(\E\big[ \| \bfA \|^\alpha\big])^{n-k}\Big(\E\big[\|\bR_{k-1}^0\|^\alpha\big]\, \mu^{(\bfA, \bfB)} \big( \big\{ (\ba, \bb) \: : \: \|\ba\| >\wt c/2\big\}\big)+ {\rm const}
\Big)\\
& \le &{\rm const} (\E\big[ \| \bfA \|^\alpha\big])^{n-k}\,,
\eeao
where we used \eqref{eq:d} in the last step. Now \eqref{eq:a} is immediate
in view of condition $\E[\|\bfA\|^\alpha]<1$ and since $f\le M$. The proof 
that $\nu$ is a Radon \ms\ on $\bbr^d_{\bf0}$ follows along the same lines.
The proof of
\beao
\lim_{n\to \infty}   \int f d \Big( \sum_{ k > n/2} \eta_{k}\Big) &=& \lim_{\nto}\sum_{k=[n/2]+1}^{ \infty} \E \Big[ \dint  
f\big( \mathbf \Pi_{k} (\ba \bR_0 + \bb)\big) \,  \mu^{(\bfA, \bfB)}(d (\ba,\bb))\Big]= 0\,,
\eeao
is an immediate con\seq\ of this fact, proving {\bf (3)}.\\[2mm]
{\bf Approximation (2).}
We have 
\beam\lefteqn{\Big| \int f d \Big( \sum_{ n/2< k \leq n} \eta_{n,k} - \sum_{ 0< k \leq n/2} \eta_{k} \Big) \Big|}\nonumber\\ &=& \Big| \int f d  \Big( \sum_{ n/2< k \leq n} (\eta_{n,k} -\eta_{n-k})\Big) \Big|\nonumber\\& =& \Big|\sum_{k= [n/2]+1}^{n}\Big(
\E  \Big[\dint  f\big( \mathbf \Pi^\downarrow_{n,k+1} 
(\ba \bR_{k-1}^0 + \bb)\big) \, \mu^{(\bfA, \bfB)}(d (\ba, d\bb))\Big]\nonumber\\
&&-\E\Big[  \dint  f\big( \mathbf \Pi_{n-k} \,(\ba \,\bR_0 + \bb)\big)\,   
\mu^{(\bfA, \bfB)}(d (\ba,\bb))\Big]\Big)\Big|\,,\label{eq:c}
\eeam
and we will show that the \rhs\ converges to zero as $\nto$.
By uniform continuity of $f$, 
\beam\label{eq:unifcont}
\mbox{for any $\vep>0$ there is $\delta>0$ such that 
$\| \bs - \br \| \leq \delta \Longrightarrow| f(\br) - f(\bs)| \leq \varepsilon$.}
\eeam
Let $(\mathbf\Pi_i')$ be an independent copy of $(\mathbf\Pi_i)$. 
For  $[n/2] < k \leq n$ write
\begin{equation*}
A_{k,\delta}( \ba) = 
\Big\{\Big \| \mathbf\Pi_{n-k}' \,\ba \,\sum_{j=k-1}^{\infty} \mathbf \Pi_j\,\bfB_{j+1} \Big\| > \delta\Big\}\,.
\end{equation*}
Since $\mathbf \Pi^\downarrow_{n,k+1} \stackrel{d}{=} \mathbf \Pi_{n-k}$ we have
\beao
\lefteqn{\Big| \E\Big[  \dint  
f\big(\mathbf \Pi^\downarrow_{n,k+1}(\ba \bR_{k-1}^0 + \bb)\big)\,   
\mu^{(\bfA, \bfB)}(d (\ba, \bb))\Big] - 
\E\Big[  \dint  f\big( \mathbf\Pi_{n-k} (\ba \bR_0 + \bb)\big)\,   
\mu^{(\bfA, \bfB)}(d (\ba, \bb))\Big]\Big|} \\
&\leq& \Big(\dint_{A_{k,\delta}(\bfa)}+\dint_{A_{k,\delta}^c(\bfa)}\Big)\\
&& \E\Big[\Big|  f\Big( \mathbf\Pi_{n-k}'\,(\ba \sum_{j=0}^{k-2} \mathbf\Pi_j\bfB_{j+1} + \bb)\Big) -  
f\Big( \mathbf\Pi_{n-k}'  (\ba\sum_{j=0}^{\infty} \mathbf\Pi_j\bfB_{j+1}+ \bb)\Big) \Big] 
\Big|  \,\mu^{(\bfA, \bfB)}(d (\ba,\bb))\\
& =& H_k^{(1)}+H_k^{(2)}.
\eeao
The following bounds hold
\beao
H_k^{(1)} &\leq & 2\,M \,\E\Big[ \mu^{(\bfA, \bfB)}  \Big(  
\Big\{ (\ba, \bb) \: : \: \| \mathbf\Pi_{n-k}'  \| \,\| \ba \| \,
\sum_{j=k-1}^{\infty} \| \mathbf\Pi_j\bfB_{j+1} \| > \delta \Big\}\Big)   \\
&				\leq & 2\,M\, 
( \E \big[\| \bfA \|^\alpha\big])^{ n-1}\, \E[ R^\alpha]\, 
\delta^{-\alpha}\,  \mu^{(\bfA, \bfB)}( \{(\ba, \bb) \: : \: \| \ba\| >1 \})\,,\\
&=&{\rm const}\,( \E \big[\| \bfA \|^\alpha\big])^{ n-1}\,\delta^{-\alpha}\,.
\eeao
Using the continuity of $f$, we also have
\beao
H_k^{(2)} &= &   \int \E\big[| \cdots  | \, \1\big(A^c_{k,\delta}(\ba)\,,
\| \mathbf \Pi_{n-k}'\|\, \|\ba\|\, R >c\big)\, \mu^{(\bfA, \bfB)}(d (\ba,\bb))\big] \\
& \leq & \varepsilon\, (\E[\|\bfA\|^\alpha])^{n-k}\,\E[ R^\alpha]\, 
\mu^{(\bfA, \bfB)}( \{ (\ba, \bb) \: : \: \| \ba\| >c\})\,.
\eeao
These computations yield
\begin{equation*}
\sum_{k= \lfloor n/2 \rfloor +1}^n \big( H_k^{(1)}+H_k^{(2)}\big) 
\leq  {\rm const}\, n  \,\big( \E \| \bfA \|^\alpha)^{ n-1} 
\delta^{-\alpha} + \varepsilon\, (\E[\|\bfA\|^\alpha])^{n-k}\big)\,.
		\end{equation*}
This bound yields that the \rhs\ of \eqref{eq:c} converges to zero by
first letting $\nto$ and then $\vep\to 0$.
\end{proof}

\begin{proof}[ Final steps in the proof of Theorem~\ref{thm:main2}]
Choose $f \in C_c^+(\R_{\mathbf 0}^d)$ and fix constants $c, M>0$ 
such that \eqref{eq:suppf} holds. 
 By uniform continuity of $f$, we can choose $\vep,\delta>0$ such that 
\eqref{eq:unifcont} holds.
Write 
\beao
A_{n,t}= \Big\{\Big\|\sum_{j=n}^{\infty} \mathbf \Pi_j\bfB_{j+1}  \Big\|>  \delta t\Big\}\,.
\eeao
We have
\beao
\big|\E\big[  f(t^{-1} \bR)  - f(t^{-1} \bR_n^0)\big]\big| &\le  &  
\E \big[ \big| f(t^{-1} \bR)  - f( t^{-1}\bR_n^0)\big|\big] \\
&=&\E\Big[ \Big| f(t^{-1} \bR)  - f\Big(t^{-1} \sum_{j=0}^{n-1} \mathbf\Pi_j\bfB_{j+1}\Big)\Big|\big(\1(A_{n,t})+ \1(A^c_{n,t})\big)\Big]\\
&= & H_1(t) + H_2(t).
\eeao
		Both terms are \asy ally negligible. Indeed, for the first one,
\beao
\dfrac{H_1(t)}{\P(\|(\bfA, \bfB)\|>t)} &\leq& 2\,M \,
\dfrac{\P(\| \mathbf \Pi_n'\|\, R >\delta t)}{\P(\|(\bfA, \bfB)\|>t)}\\
 &\leq &{\rm const}\, (\E \|\bfA \|^\alpha)^n\,\delta^{-\alpha} \,.
\eeao
The \rhs\ converges to zero by first letting $t\to\infty$ and then $\nto$,
also observing that $\E[\|\bfA\|^\alpha]<1$. For the second one,  using 
\eqref{eq:unifcont},
\beao
\dfrac{H_2(t)}{\P(\|(\bfA, \bfB)\|>t)} & = &  \dfrac{
\E\Big[ \Big| f(t^{-1} \bR)  - 
f\Big(t^{-1} \sum_{j=0}^{n-1} \Pi_j\bfB_{j+1}\Big)\Big|\1\big(A^c_{n,t})\,\1(\| \bfR \|>c\,t)\Big]}{\P(\|(\bfA, \bfB)\|>t)} \\
&\leq & \varepsilon\, \dfrac{ \P( \| \bfR \| > ct)}{\P(\|(\bfA, \bfB)\|>t)}\le {\rm const}\,\vep\,.
\eeao
		In view of Lemma~\ref{lem:last} we may conclude that if we first take $t\to \infty$, then $n \to \infty$ followed by $\varepsilon \to 0$, we may conclude that
		\begin{equation*}
			\frac{\E  f(t^{-1} \bR) }{\P(\|(\bfA, \bfB)\| >t)} \to \int f(\br) \nu(d\br).
		\end{equation*} 
		Since $f$ is arbitrary the theorem follows.
	\end{proof}

	\appendix
\section{}
\subsection{Proof of Lemma~\ref{lem:productregvaruniv}}\label{app:prooflem1.2}
(1) was proved in Embrechts and Goldie \cite{embrechts:goldie:1980}, p.~245.
	We start with (2). Observe that for any $M>0$, by the uniform
	\con\ theorem for \regvary\ \fct s,
	\beao
	\dfrac{\P(XY>x)}{\P(X>x)}\ge \int_0^M \dfrac{\P(X>x/y)}{\P(X>x)}\,\P(Y\in dy)
	\to \int_0^M y^\alpha\,\P(Y\in dy)\,,\qquad \xto\,.
	\eeao
	If $\E[Y^\alpha]=\infty$ we can make the \rhs\ arbitrarily large by letting 
	$M\to\infty$.
	\par
	We continue with (3). We follow the lines of the proof 
	of Proposition~3.1 in Davis and Resnick \cite{davis:resnick:1985a}
	who consider the case of iid $X,Y$. Choose any $M>1$. Then
	\beao
	&&\P(XY>t)\\&=& \P(XY>t\,,X\le M)+ \P(XY>t\,,M<X\le t/M)+ \P(XY>t\,,X>t/M)\\
	&=&I_1(t)+I_2(t)+I_3(t)\,.
	\eeao
	In view of \eqref{eq:strange}, $I_2(t)/\P(X>t)$ is \asy ally negligible
	when first $t\to\infty$ and then $M\to\infty$. In view of Breiman's
	Lemma~\ref{lem:breiman} we have as $t\to\infty$,
	\beao
	\dfrac{I_3(t)}{\P(X>t)}= \frac{\P( X (Y \wedge M)>t)}{\P(X>t)}  &\to& \E\big[(Y\wedge M)^\alpha]\,, \\
	\dfrac{I_1(t)}{\P(X>t)}= \dfrac{\P\big(Y\,X \1( X \leq M)>t \big)}{\P(Y>t)}\dfrac{\P(Y>t)}{\P(X>t)}
	&\to & c_0\,\E\big[X^\alpha\1_{\{ X \leq M \}}\big]\,,\\
	\eeao
	where $c_0=\lim_{t\to\infty}\P(Y>t)/\P(X>t)$ is assumed finite. Now the desired
	result follows when $M\to\infty$.
\subsection{A result from \cite{damek:dyszewski:2018}}
\ble\label{lem:dd2018}
Assume that $\| \bfA\|$ is \regvary\ with index $\alpha>0$, 
$\E[\|\bfA\|^\alpha] <1$, $\P(\| \bfB\|>t) = O(\P(\|\bfA\|>t))$,
and
\begin{equation*}
\dfrac{\P(\|\bfA_1\| \cdot \|\bfA_2\|>t)}{\P(\|\bfA\|>t)} \to 2\, \E [\|\bfA\|^\alpha]\,,\qquad t\to\infty\,.
\end{equation*}
Then $R = \sum_{k=0}^\infty \| \bfB_{k+1}\|\, \prod_{j=1}^k \| \bfA_j\|$ is
finite and satisfies $\P(R>t) = O(\P(\|\bfA\|>t))$ as $t\to\infty$. 
In particular, $\E [R^\alpha]<\infty$.
\ele


\end{document}